\documentstyle{article}

\newtheorem{theorem}{Theorem}[section]
\newtheorem{lemma}[theorem]{Lemma}
\newtheorem{prop}[theorem]{Proposition}
\newtheorem{cor}[theorem]{Corollary}

{\catcode`\@=11
\gdef\setft#1#2#3e{%
\def\@oddfoot{
{\setbox0=\hbox{#1}
\setbox1=\hbox{#3}
\ifdim\wd0>\wd1
\dimen0=\wd0
\box0\hfil#2\hfil\hbox to\dimen0{\hfil\hfil\box1}
\else \dimen0=\wd1
\hbox to\dimen0{\box0\hfil }\hfil#2\hfil\box1 \fi
}}} }

\title{On Tangents to Quadric Surfaces}
\author{Ciprian Borcea, Xavier Goaoc, Sylvain Lazard and Sylvain Petitjean}
\date{}

\begin{document}
\maketitle

\noindent
\begin{abstract}
\noindent
We study the variety of common tangents for up to four quadric surfaces
in projective three-space, with particular regard to configurations of
four quadrics admitting a continuum of common tangents.

\medskip \noindent
We formulate geometrical conditions in the projective space defined by
all complex quadric surfaces which express the fact that several quadrics
are tangent along a curve to one and the same quadric of rank at
least three, and called, for intuitive reasons: {\em a basket}.
Lines in any ruling of the latter will be common tangents.

\medskip \noindent
These considerations are then restricted to spheres in Euclidean three-space,
and result in a complete answer to the question {\em over the reals}:
``When do four spheres allow infinitely many common tangents?''.

\end{abstract}


\normalsize

\section*{Introduction}

\medskip \noindent
Tangents to a non-singular complex projective quadric surface make-up a
threefold, namely: the projectivized tangent bundle of the given quadric.
After a birational contraction, this threefold can be represented as a
{\em quadratic section} of the Grassmannian $G(2,4)$ of all projective
lines in $P_3(C)$ i.e. all 2-subspaces of the vector space $C^4$.

\medskip \noindent
$G(2,4)$, in its Pl\"{u}cker embedding, is itself a quadric in $P_5(C)$,
and it follows that {\em four non-singular quadric surfaces in general
position allow $2^5=32$ common tangents}.

\medskip \noindent
However, what we may call {\em degenerate configurations} of four quadrics,
would still allow a continuum of common tangents. This obviously happens
when the four quadrics have a common curve of intersection, other than
a union of less than four lines (or, as we shall observe later, when their
{\em duals} do).
In general, a curve of common tangents would give a {\em ruled surface} in
$P_3(C)$, which is tangent along some curve with each of the given quadric
surfaces.

\medskip \noindent
{\em We investigate the case when this ruled surface is itself a quadric
surface}. A simple example (where one can ``see'' what happens
for the {\em real points}) is that of a hyperboloid of revolution
in which one throws four spherical ``balls''  and lets them
rest when reaching a circle of tangency with their ``basket''.
Because of this intuitive background, we introduce:

\medskip \noindent
{\bf Definition:}\
Let $q_1$ and $q_2$ be distinct quadric surfaces of rank at least three
(i.e. non-singular or with at most an isolated conic singularity). We say
that $q_1$ is a {\bf basket} for $q_2$ (and then, $q_2$ will be a
basket for $q_1$) when the two quadrics are tangent along a conic.
(Thus, the intersection of the two quadrics is represented by twice this
conic.)

\medskip \noindent
{\bf Convention:}\ {\em In the sequel, whenever we speak of a common
basket $b$ for quadrics $q_i$, we assume that all quadrics concerned
are distinct and of rank at least three.}

\medskip \noindent
We are going to formulate conditions expressing the fact that two, three
or four quadrics allow a common basket. These will be geometrical
conditions in the space of all quadric surfaces which is a nine-dimensional
complex projective space corresponding to lines through zero in the
vector space of all $4\times 4$ symmetric matrices with complex entries:

$$ P_9(C)=P(Sym_C(4)) $$

\noindent
Indeed, one may identify quadratic forms and symmetric matrices (over C) via
the standard bilinear form $<,>$:

$$ q(x)=\sum_{i,j} q_{ij}x_ix_j =<x,Qx> \ \ \ \ \ \ Q=Q^t $$

\noindent
The rank of the quadric $q$, as already spoken of, is simply the rank of
the matrix $Q$.

\medskip \noindent
Sometimes, we'll refer to quadrics of rank at most three as {\em cones},
while those of rank at most two, respectively one, will be called
{\em two-planes}, respectively {\em double-planes}. Obviously, a two-plane
is a cone over a degenerate conic i.e. a two-line.

\medskip \noindent
The closure of the rank three locus is denoted ${\cal R}^3_8$,
the closure of the rank two locus is denoted ${\cal R}^2_6$,
and the rank one locus is denoted ${\cal R}^1_3$.

\medskip \noindent
To begin with, we give an equivalent of our definition of baskets:

\begin{prop}
$b$ is a basket for $q$ if and only if the pencil
$[b,q]=\lambda b + \mu q$
meets the rank one locus (i.e. contains a double-plane).
\end{prop}

\medskip \noindent
{\em Proof:}\ If the pencil contains a double plane, the two quadrics
intersect in a double conic and must be tangent along it.

\medskip \noindent
Conversely, the double plane through the conic of tangency has the
same intersection with $b$ as $q$ (and with $q$ as $b$), and must
belong to the pencil. \ \ \ $\Box$

\medskip \noindent
In the same vein, we shall obtain:

\begin{prop}
The quasi-projective variety (see our convention above):

$$ B^2_{16}=\{ (q_1,q_2,b)\ :\ b \ \mbox{is a common basket for} \ q_1 \
\mbox{and}\  q_2 \  \} \subset (P_9-{\cal R}^2_6)^3 $$

\noindent is irreducible, of dimension sixteen, and the pair $(q_1,q_2)$ of a generic point $(q_1,q_2,b)$ is
characterized by the property that the pencil $[q_1,q_2]$ contains a two-plane.
\end{prop}

\begin{prop}
The quasi-projective variety:

$$ B^3_{21}=\{ (q_1,q_2,q_3,b)\ :\ b\ \mbox{is a common basket for} \
q_1, q_2 \ \mbox{and}\  q_3 \  \} \subset (P_9-{\cal R}^2_6)^4 $$

\noindent is irreducible, of dimension twenty one, and the triple $(q_1,q_2,q_3)$ of a generic point
$(q_1,q_2,q_3,b)$ is characterized by the property that the span $[q_1,q_2,q_3]\approx P_2$ contains a pencil of
cones with the same vertex, and the rank two points in this pencil are precisely where it meets the lines
$[q_i,q_j]$.
\end{prop}

\begin{prop}
The quasi-projective variety:

$$ B^4_{25}=\{ (q_1,q_2,q_3,q_4,b) \ : \ b\ \mbox{is a common basket for}\
q_i, \ i=1,...,4 \  \} \subset (P_9-{\cal R}^2_6)^5 $$

\noindent
is irreducible, of dimension twenty five.

\medskip \noindent
The quadruple $(q_1,...,q_4)$ of a generic point $(q_1,...,q_4,b)\in B^4_{25}$ is characterized by the following
property:\ the span $[q_1,..,q_4]\approx P_3$ contains a complete quadrilateral consisting of four pencils of
cones; the six vertices $p_{kl}$ lie, respectively, on the six lines $[q_i,q_j]$ and correspond precisely with
the rank two quadrics of the quadrilateral.
\end{prop}

\medskip \noindent
When we present our proofs, we'll examine and characterize all
possibilities {\em in terms of conditions on the configuration} $(q_i)$,
indicating how to `reconstruct' all baskets when the conditions are met.

\medskip \noindent
These results then lead to a proof of {\em uniqueness}, up to the action
of the projective automorphism group $PSL_C(4)$, of a ``double-four''
configuration, namely: two {\em quadruples of linearly independent and
smooth quadrics}, such that each quadric in one group is a common
basket for the quadrics in the other group. This configuration arises
quite naturally from the point of view of duality, and is related to
the {\em Reye configuration} $(12_4,16_3)$.

\medskip \noindent
When seen on the Grassmannian $G(2,4)$, the tangents of a smooth quadric
surface define a {\em degenerate quadratic line complex}: it has
singularities along the two conics representing the two rulings of the
quadric.

\medskip \noindent
The relation of generic quadratic line complexes with Kummer
surfaces is a classical subject \cite{Jes}, \cite{Hud}, \cite{GH}.
In our case, Kummer surfaces appear when intersecting two degenerate
quadratic line complexes, that is, when considering common tangents
to a generic pair of quadrics. A generalisation to higher-dimensional
Calabi-Yau varieties is pursued in \cite{Bor2}.

\medskip \noindent
In the last section we use our results on quadrics with common baskets
to solve the problem of describing {\em all possible configurations of
four spheres in $R^3$ with infinitely many real common tangents}. The
conclusion agrees with intuitive expectations: {\em the four centers have to
be collinear}, and the radii must accommodate one of the following
possibilities:

\medskip
(i)\ \ \  the four spheres intersect in a common circle or point;

(ii)\ \ the radii are equal, and there's a common cylindrical basket;

(iii)\ the four spheres have a common conical basket;

(iv)\ \ there's a common basket in the shape of a hyperboloid of revolution
with one sheet and axis the line of centers.

\medskip \noindent
This complements results in \cite{MPT} and \cite{ST} on configurations
of four spheres with a finite number of common tangents. The effective
upper bound is 12.

\medskip \noindent
The material is organized in eight sections:

\medskip
1.\ Stratification by rank in $P_9=P(Sym_C(4))$

2.\ Two quadrics in a basket

3.\ Three quadrics in a basket

4.\ Four quadrics in a basket

5.\ A double-four example

6.\ Tangents and Grassmannians

7.\ Duality

8.\ Common tangents to four spheres in $R^3$

\section{Stratification by rank in $P_9=P(Sym_C(4))$}

\medskip \noindent
In this section we review some classical facts about the space of
all quadric surfaces. More general considerations can be found in
\cite{Bor1}, or \cite{Har}.

\medskip \noindent
The stratification by rank yields, in the case of $4\times 4$
symmetric matrices three determinantal varieties:

$$ {\cal R}^1_3\subset {\cal R}^2_6\subset {\cal R}^3_8\subset
P_9=P(Sym_C(4)) $$

\noindent
The {\em rank at most three locus} ${\cal R}^3_8$ is the {\em degree
four hypersurface} defined by all singular quadrics:

$$ {\cal R}^3_8=\{ Q\in P_9\ : \ det(Q)=0 \} $$

\noindent
These singular quadrics are obviously {\em cones} over conics in some
plane $P_2\subset P_3$ constructed from a vertex outside that plane.
Generically, the vertex is the only singularity.

\medskip \noindent
The {\em rank at most two locus} ${\cal R}^2_6$ is codimension two
in ${\cal R}^3_8$ and represents in fact its singular locus. It is
defined in $P_9$ by the vanishing of all $3\times 3$ minors, and
Giambelli's formula gives its degree as ten.

\medskip \noindent
A quadric in ${\cal R}^2_6$ is the cone over some degenerate conic,
that is: two lines, and so the union of {\em two planes}.

\medskip \noindent
The {\em rank one locus} ${\cal R}^1_3$ is codimension three in
${\cal R}^2_6$ and represents its singular locus. It is defined in
$P_9$ by the vanishing of all $2\times 2$ minors, and can also
be described as the {\em image of the quadratic Veronese embedding}:

$$ v: P_3 \rightarrow P_9, \ \ \ v(x)= x^t\cdot x $$

\noindent
where $x=(x_0:x_1:x_2:x_3)$, and $x^t$ stands for its (column) transpose.
Clearly, the symmetric matrix $x^t\cdot x$ has rank one, and:

$$ {\cal R}^1_3= v(P_3) \subset P_9 $$

\medskip \noindent
Similarly, ${\cal R}^2_6$ can be identified with the quotient of
$P_3\times P_3$ by the involution $\sigma(x,y)=(y,x)$ using:

$$ w: (P_3)^2/\sigma \rightarrow P_9, \ \ \ w(x,y)=w(y,x)=
\frac{1}{2}(x^t\cdot y + y^t\cdot x) $$

\noindent This shows that ${\cal R}^2_6$ is swept out by a three parameter (rational) family of projective
three-spaces in $P_9$. It is also swept out by a family of projective two-spaces in $P_9$, indexed by the
Grassmannian $G(2,4)$: indeed, given $\ell\in G(2,4)$, which we regard as a two-subspace of $C^4$, we have a
plane in ${\cal R}^2_6$ given by:

$$ P_2\approx \{ Q \ : \ \ell\subset ker(Q)\} \subset {\cal R}^2_6 $$

\medskip \noindent
Later considerations will involve various pencils of singular quadrics,
and it may be remarked here that the first of the above families
provides a seven parameter family of pencils (i.e. lines) in
${\cal R}^2_6$ , to be denoted by ${\cal M}^2_7\subset G(2,10)$,
while the second family provides a six parameter family of pencils
in ${\cal R}^2_6$, to be denoted by ${\cal F}^2_6\subset G(2,10)$.

\medskip \noindent
Here, $G(2,10)$ stands for the Grassmannian of all lines in $P_9$,
where we can ask for the intersection of ${\cal M}^2_7$ and
${\cal F}^2_6$. In terms of symmetric matrices, pencils in the
first family involve quadrics with a given vector in the image,
while pencils in the second family involve quadrics with a given
two-space in the kernel. Considering that $im(Q)=ker(Q)^{\perp}$,
we see that the intersection is five dimensional:

$$ {\cal M}^2_7 \cap {\cal F}^2_6 = {\cal T}_5 \subset G(2,10) $$

\noindent and consists of pencils of quadrics parametrized by pairs $(\ell,x)\in G(2,4)\times P_3$, with
$x\subset \ell^{\perp}$, and defined by:

$$ P_1\approx \{ Q \ : \ \ell\subset ker(Q), x\subset im(Q) \}
\subset {\cal R}^2_6 $$

\medskip \noindent
These relations can be observed from another point of view, based
on the fact that ${\cal R}^2_6$ is precisely the {\em secant variety}
of ${\cal R}^1_3\subset P_9$. Indeed, ${\cal R}^2_6$ can be obtained
as the closure of the union of all lines spanned by two distinct
points in ${\cal R}^1_3$. The closure brings in the points of all
lines tangent to ${\cal R}^1_3$, and one can see that the family
${\cal F}^2_6$ consists of secants and tangents to ${\cal R}^1_3$,
while ${\cal T}_5\subset {\cal F}^2_6$ retains only the lines tangent
to the rank one locus.

\medskip \noindent
We note that an arbitrary line in $P_9$ can meet the rank one locus
in at most two points, since any secant is, in adequate coordinates,
of the form: \ $\lambda_1x_1^2+\lambda_2x_2^2$, and this pencil has
no other double-plane.

\medskip \noindent
The family of projective three-spaces sweeping out ${\cal R}^2_6$ can be
recognized now as the family of all tangent spaces to the rank one locus,
and this makes obvious our earlier result:

$$ {\cal M}^2_7 \cap {\cal F}^2_6 = {\cal T}_5 \subset G(2,10) $$

\noindent
since a line in a tangent space to ${\cal R}^1_3$ must pass through
the point of tangency in order to be in ${\cal F}^2_6$, and is then
perforce in ${\cal T}_5$. \ \ \ $\Box$

\medskip \noindent
We can expand our description of pencils of singular quadrics by
considering those contained in ${\cal R}^3_8$ and not already
contained in ${\cal R}^2_6$. The generic quadric in such a pencil
will have a unique singular point (the vertex of the cone), and
we may expect two types of pencils:

\medskip
(F)\  with {\em fixed} vertex,\ \  or \ \ (M)\ with {\em moving} vertex.

\medskip \noindent
The first type obviously arises by choosing a pencil of conics in
some $P_2\subset P_3$ and constructing the cones over the conics
in the pencil from a fixed vertex away from our $P_2$. This yields
an eleven parameter family, to be denoted ${\cal F}^3_{11}$.

\medskip \noindent
This family is related to the fact that ${\cal R}^3_8$ is swept out
by a family of projective five-spaces indexed by $P_3$. Indeed, for
$x\in P_3$, we have:

$$ P_5\approx \{ Q\ :\ Qx=0\} $$

\noindent
which describes all quadrics singular at $x$. Lines in this $P_5$
make-up a Grassmannian $G(2,6)$ of dimension eight, hence our eleven
dimensional ${\cal F}^3_{11}$.

\medskip \noindent
Yet another description of this family is related to the fact that
${\cal R}^3_8$ is the {\em variety of secant planes} of
${\cal R}^1_3\subset P_9$.
Indeed, ${\cal R}^3_8$ is the closure of the union of all
planes spanned by three double-planes. Thus, ${\cal R}^3_8$ is
also swept out by a nine parameter family of planes, hence our
${\cal F}^3_{11}$, made of lines in these planes.

\medskip \noindent
For the second type (M), let us observe first that if we write the
pencil with moving vertex as:\ $\lambda_1q_1+\lambda_2q_2$, with
$q_i$ of rank three and with vertex $v_i$, then all cones in the
pencil must contain the line $[v_1,v_2]$. This follows from the fact
that the tangent hyperplane to ${\cal R}^3_8$ at $q_i$ consists
of all quadrics passing through $v_i$, and our pencil lies in the
intersection of these two hyperplanes. Thus, $q_i$ passes through
$v_j$, and contains $[v_i,v_j]$.

\medskip \noindent
Clearly, the vertices in the pencil must move along this line,
and the intersection $q_1\cap q_2$ will consist of a conic and
the double line $[v_1,v_2]$.

\medskip \noindent
Thus, a pencil of type (M) arises by considering a fixed conic in some $P_2\subset P_3$, then choosing a line
through a point of the conic, but not contained in its plane, and moving a vertex (linearly) along this line.
This yields an eleven parameter family, to be denoted ${\cal M}^3_{11}$.

\medskip \noindent
We summarize these results in:

\begin{prop}
The variety of lines contained in ${\cal R}^3_8\subset P_9$ consists of two irreducible components of dimension
eleven: ${\cal F}^3_{11}$ and  ${\cal M}^3_{11}$, made generically of pencils with fixed singularity,
respectively moving singularity.

\medskip \noindent
The variety of lines contained in ${\cal R}^2_6\subset P_9$ consists
of two irreducible components: ${\cal F}^2_6$ of dimension six, and
${\cal M}^2_7$ of dimension seven.

\medskip \noindent
We have inclusions:

$$ {\cal F}^2_6\subset {\cal F}^3_{11}\ \ \mbox{and} \ \
{\cal M}^2_7\subset  {\cal M}^3_{11} $$

\noindent
and intersections:

$$ {\cal F}^3_{11} \cap {\cal M}^3_{11} =
{\cal F}^2_6 \cap {\cal M}^2_7 = {\cal T}_5 $$

\noindent
with ${\cal T}_5\subset G(2,10)$ standing for the five dimensional
variety made of tangent lines to the rank one locus ${\cal R}^1_3$.
\ \ \ $\Box$
\end{prop}

\medskip \noindent
{\bf Remark:}\ We have the following implications:

\medskip
if $\ell\in {\cal M}^3_{11}-{\cal M}^2_7$, then $\ell$ has a single two-plane;

\medskip
if $\ell\in {\cal M}^2_7-{\cal T}_5$, then $\ell$ has no double-plane;

\medskip
if $\ell\in {\cal F}^3_{11}-{\cal F}^2_6$, then $\ell$ has three two-planes, counting multiplicity;

\medskip
if $\ell\in {\cal F}^2_6-{\cal T}_5$, then $\ell$ has two double-planes.

\section{Two quadrics in a basket}

\medskip \noindent
In this section we prove and refine Proposition 0.2.

\medskip \noindent
Suppose $q_1$ and $q_2$ allow a common basket $b$. Then, by
Proposition 0.1, the pencil $[b,q_i]$ meets the rank one locus
in $d_i$. Then, either $b,q_1,q_2$ are collinear, and so $[q_1,q_2]$
meets the rank one locus in $d=d_1=d_2$, or $[q_1,q_2]$ meets the
rank two locus where it intersects $[d_1,d_2]$.

\medskip \noindent
Conversely, if $[q_1,q_2]$ meets the rank one locus, then any $b$
(of rank at least three) in the pencil will be a common basket.
If $[q_1,q_2]$ meets the rank two locus in $p$, then there's a secant
of the rank one locus through $p$ and two double-planes we label $d_i$.
Then $[q_i,d_i], \ i=1,2$ are coplanar and meet in a point $b$ which
is a common basket for $q_i$.

\medskip \noindent
Irreducibility and the dimension count for $B^2_{16}$ follows from
the fact that there's an eight parameter family of pencils through
each point of ${\cal R}^2_6$, and on each pencil, the choice of two
points means two more parameters

\medskip \noindent
This already proves Proposition 0.2, but we may refine the statement by observing the `reconstruction' process
of a common basket in more detail.

\medskip \noindent
As a {\bf rule}, whenever the variety of common baskets has positive
dimension, we'll describe its {\em closure}, being understood that,
according to our convention, we retain only the generic part made
of quadrics of rank at least three for the role of baskets.

\medskip \noindent
Thus, when $[q_1,q_2]$ meets the rank one locus, the whole pencil
offers common baskets, but there is one type of situation where
we have in addition, another rational curve of common baskets:
let us call $d$ the double-plane on $[q_1,q_2]$, and suppose that
there's a common basket $b$ away from this pencil. Then, we get
double-planes $d_i$ on $[b,q_i]$, and $[q_1,q_2]$ is in the
plane $[d,d_1,d_2]$ which lies in ${\cal R}^3_8$.

\medskip \noindent
It follows that, in our situation, $[q_1,q_2]$ is a pencil of
cones with fixed vertex, and as $q_i$ is also a basket for $q_j$,
the two quadrics are cones from the same vertex over two
(non-singular) conics which have two points of tangency.
We note that $[q_1,q_2]$ has a rank two point at the intersection
with $[d_1,d_2]$, and we may run other secants of the rank one
locus through this point and construct other common baskets.

\medskip \noindent
The fact that this leads to another rational curve of common
baskets follows from the same argument as in the reconstruction
process for the hypothesis of $[q_1,q_2]$ meeting the rank two
locus, considered presently.

\medskip \noindent
We'll need some lemmas:

\begin{lemma}
If a plane $P_2\subset P_9=P(Sym_C(4))$ contains four distinct
rank one quadrics, then the plane contains a (non-singular) conic
of rank one quadrics. More precisely, $P_2$ is then the span of
the Veronese  image $v(P_1)\subset P_2$ of some line $P_1\subset P_3$.
\end{lemma}

\medskip \noindent
{\em Proof:}\ When we look at the planes in $P_3$ corresponding
to our four double-planes, we see that they must have a point
in common, otherwise they would be projectively equivalent to
$x_i^2=0, \ i=0,1,2,3$, and the double-planes would span a
$P_3\subset P_9$. Thus, the problem is reduced to its version
in $P_5=P(Sym_C(3))$.

\medskip \noindent
Now, again, the four double-lines must have a common point
in $P_2$, for otherwise three of them would be projectively
equivalent to $x_i^2=0,\ i=0,1,2$ and their span in $P_5$
has no other double-line.

\medskip \noindent
Thus, the original four double-planes have a common line and are projectively equivalent to four points in the
family $(\lambda_0x_0+\lambda_1x_1)^2, \ (\lambda_0:\lambda_1)\in P_1$ which is the Veronese image of a line. \
\ \ $\Box$

\begin{lemma}
Let $p\in {\cal R}^2_6-{\cal R}^1_3$ be a two-plane. The family of
secants (and tangents) to the rank one locus which pass through $p$
make up a rational curve.

\medskip \noindent
In fact, there's a unique plane  $P_2\subset {\cal R}^2_6$ passing
through $p$ and containing (as a non-singular conic) the Veronese image
$v(P_1)\subset P_2$ of a line. Thus, all lines through $p$ in this
$P_2$ make up a rational curve of secants of the rank one locus,
with exactly two tangents. \ \ \ $\Box$
\end{lemma}

\medskip \noindent
Now, whenever we have a rank two point $p\in [q_1,q_2]$, we can
run all proper secants through $p$ and label, in two ways, the
two rank one points on it $d_1$ and $d_2$. Each labelling gives a common
basket $b=[q_1,d_1]\cap [q_2,d_2]$. Thus, the curve of common baskets is
a double covering of $P_1$ ramified over the two tangents, and
hence a rational curve itself. The association $b\mapsto d_i=[b,q_i]\cap
v(P_1)$ gives an isomorphism with $v(P_1)$, for $i=1,2$.

\medskip \noindent
Actually, the curve of common baskets is the residual intersection
of the two cones over $v(P_1)$ with vertex at $q_1$, respectively
$q_2$, and consequently a {\em conic} itself.

\medskip \noindent
Thereby we obtain this complement to Proposition 0.2:

\begin{prop}
The variety of common baskets for two quadrics $q_i$ is
determined by the number of points in the intersection
$[q_1,q_2]\cap {\cal R}^2_6$, and consists of as many
rational curves.

\medskip \noindent
In order to have more than one rational curve, it is necessary
and sufficient that $[q_1,q_2]$ be a pencil in
${\cal F}^3_{11}-{\cal F}^2_6$, that is: $q_1$ and $q_2$ must be two cones
with the same vertex, and in this case we have (counting multiplicities)
three rational components. \ \ \ $\Box$
\end{prop}

\section{Three quadrics in a basket}

\medskip \noindent
In this section we prove and expand Proposition 0.3.

\begin{lemma}
Let $\ell$ be a pencil of cones with three distinct rank two points, but not contained in ${\cal R}^2_6$. Then
there's a unique trio of double-planes $d_i, \ i=1,2,3$ (up to permutation), such that the span $P_2\approx
[d_1,d_2,d_3]$ contains the pencil.
\end{lemma}

\medskip \noindent
{\em Proof:}\ From section 1 it follows that $\ell$ consists of cones with a fixed vertex over a pencil of
conics with three distinct two-lines. All such pencils of conics are equivalent under projective transformations
of the plane (i.e. under $PSL_C(3)$), and thus the pencil $\ell$ can be turned into the diagonal form: \
$\lambda(x_1^2+x_2^2)+\mu(x_2^2+x_3^2)$, which is clearly in the span of the three double-planes $d_i=x_i^2$.
This proves the existence part.

\medskip \noindent
For uniqueness (up to permutation), suppose we have another trio of double-planes $d'_i$ with $\ell\subset
[d'_1,d'_2,d'_3]$. Then $\ell=[d_1,d_2,d_3]\cap [d'_1,d'_2,d'_3]$ and with proper indexing, the edges
$[d_i,d_j]$ and $[d'_i,d'_j]$ meet $\ell$ in the same point of rank two $p_{ij}$.

\medskip \noindent
By the reciprocal of Desargues' theorem (in dimension three)
the triangles determined by $d_i$, respectively $d'_i$, are in
perspective i.e. the lines $[d_i,d'_i]$ meet at a point $p$,
which is necessarily of rank two. But Lemma 2.1 requires then
all our double-planes to be on the same conic. This contradiction
proves the uniqueness part. \ \ \ $\Box$

\medskip \noindent
Let $q_i, \ i=1,2,3$ be three distinct quadrics with a
common basket $b$. Again, the pencils $[b,q_i]$ must meet
the rank one locus in double-planes $d_i$.

\medskip \noindent
Suppose first that $q_i$ are collinear. Then, either $b$ is on the
same line and $d_1=d_2=d_3$, or $b$ is away from this line and then
$d_i$ are distinct and span a plane containing the line.

\medskip \noindent
The latter case means that $q_i$ belong to a pencil $\ell$ of cones with fixed vertex which has its rank two
points at the intersections $p_i=\ell\cap [d_j,d_k]$. Obviously the six points $p_i,q_j$ on $\ell$ must satisfy
a relation, since $q_j$ are projections of $d_j$ from $b$.

\medskip \noindent
One can guess this relation from the fact that it comes from a (rational) map $P_2 \cdot \cdot \rightarrow
(P_1)^3$ whose image should be a surface of multi-degree $(1,1,1)$. The formula should also have permutation
invariance. Indeed, the relation can be written as:

$$ (p_1,p_2;p_3,q_1)+(p_2,p_3;p_1,q_2)+(p_3,p_1;p_2,q_3)=
\frac{3}{2} \ \ \ \ \ \ \ (c1) $$

\noindent where $(a,b;u,v)$ denotes the cross-ratio of four points on a projective line (in our case
$\ell\approx P_1$):

$$ (a,b;u,v)=\frac{a-u}{a-v}\cdot \frac{b-v}{b-u} $$

\noindent
Conversely, when condition $(c1)$ is satisfied, the lines $[q_i,d_i]$
are concurrent, yielding the basket $b$.

\medskip \noindent
This settles the collinear case, and we may assume now that
the three quadrics span a plane $P_2\approx [q_1,q_2,q_3]$.

\medskip \noindent
If $b$ is on one of the edges, say $[q_i,q_j]$, we have a
double-point $d=d_i=d_j$ on this line and $d_k$ on $[b,q_k]$.
Thus, {\em the existence of the proper secant (of the rank one locus)
$[d,d_k]\subset [d_1,d_2,d_3]$ intersecting an edge in a double-plane
characterizes this case}.

\medskip \noindent
Henceforth, we shall assume that $b$ is not on the lines $[q_i,q_j]$.
Then, the double-planes $d_i$ span a plane $[d_1,d_2,d_3]$,
and either it coincides with $[q_1,q_2,q_3]$ (and contains $b$), or
the two planes meet in a line.

\medskip \noindent
Let's take first the case $b\in [q_1,q_2,q_3]=[d_1,d_2,d_3]$.
Thus $[q_1,q_2,q_3]$ is a generic secant plane of the rank one locus
(i.e. not contained in ${\cal R}^2_6$), and we look at the situation
where a common basket lies in this same plane, but away from the edges
$[q_i,q_j]$.

\medskip \noindent
Then all our quadrics  are cones with the same vertex, namely the intersection of the three planes in $P_3$
corresponding to the double-planes $d_i$. (The intersection cannot be a line since then $[d_1,d_2,d_3]$ would
contain a conic $v(P_1)$ of double-planes, and this would contradict the assumption that $q_i$ are of rank at
least three.) Since the triangles $\bigtriangleup(q_i)$ and $\bigtriangleup(d_i)$ are in perspective, Desargues'
theorem says that the corresponding edges meet in collinear points $p_{ij}=[q_i,q_j]\cap [d_i,d_j] \in \ell$.

\medskip \noindent
Thus, $[q_1,q_2,q_3]$ contains a line $\ell$ of singular quadrics, meeting the pencils $[q_i,q_j]$ in three
distinct rank two points $p_{ij}$. In view of Lemma 3.1 and by the reciprocal version of Desargues' theorem,
this is enough to ensure that the triangle of our three quadrics is in perspective with the triangle (properly
labeled) of the three double-planes, and the basket $b$ is retrieved as the point of perspective.

\medskip \noindent
Note that $[q_i,q_j]$ have themselves three rank two
points on them, hence the case under consideration arises only
when {\em there's one more collinearity amongst the nine points
$[q_i,q_j]\cap [d_k,d_l]$ besides the six edges}.

\medskip \noindent
To conclude, we take up the case of a line intersection $\ell= [q_1,q_2,q_3]\cap [d_1,d_2,d_3]$.

\medskip \noindent
Again, the points $p_{ij}=[q_i,q_j]\cap [d_i,d_j]$ are rank two points on $\ell$, and we have a
three-dimensional Desargues configuration.

\medskip \noindent
If $\ell$ is not contained in ${\cal R}^2_6$, and this is obviously the {\em generic} case envisaged in our
Proposition 0.3, then Lemma 3.1 and the fact that Desargues' theorem works in both directions (from $b$ to the
$p_{ij}$, and from $p_{ij}$ to $b$), yield the result that {\em the existence of a pencil of cones $\ell\subset
[q_1,q_2,q_3]$ with $\ell\cap [q_i,q_j]$ of rank two and the rest of rank three, characterizes this situation
(and the associated common basket is uniquely determined by $\ell$)}.

\medskip \noindent
We are left with the degenerate case where $\ell=[q_1,q_2,q_3]\cap [d_1,d_2,d_3]$ is contained in ${\cal
R}^2_6$, that is: $\ell\in {\cal F}^2_6$. In other words, $\ell$ is a secant (or tangent) to the rank one locus.

\medskip \noindent
Under our assumptions $d_i$ are not on $\ell$, and Lemma 2.1 implies that $[d_1,d_2,d_3]$ is the unique $P_2$
which is the span of a Veronese curve $v(P_1)$ and contains $\ell$.

\medskip \noindent
The question is whether this pencil of two-planes (with fixed singularity line) $\ell\subset [q_1,q_2,q_3]$,
together with its three marked points $p_{ij}=\ell\cap [q_i,q_j]$ of rank two, is sufficient information for
finding a common basket.

\medskip \noindent
The answer is given in the following lemma, which yields, counting multiplicity, two common baskets
corresponding to $(\ell, p_{ij})$:

\begin{lemma}
Let $v(P_1)\subset P_2$ be a Veronese conic of double-planes, and $\ell\subset P_2$ a line with three distinct
marked points $p_{ij}$ away from the intersection $\ell\cap v(P_1)$.

\medskip \noindent
If $\ell$ is a proper secant of the Veronese conic, there are exactly two triangles $\bigtriangleup(d_i)$ and
$\bigtriangleup(d'_i)$ with vertices on the conic, and such that their edges $[d_i,d_j]$, respectively
$[d'_i,d'_j]$, meet $\ell$ in $p_{ij}$.

\medskip \noindent
If $\ell$ is tangent to the Veronese conic, there's only one such triangle.
\end{lemma}

\medskip \noindent
{\bf Remark:} \ This is clearly related to Desargues' theorem,
but requesting the two triangles in perspective to have their vertices
on a conic. We have therefore a five parameter family with a (rational) map
to lines in $P_2$ with three marked points: another five parameter family.
One can fairly expect the map to be a birational equivalence, which
indeed turns out to be the case.

\medskip \noindent
There's an alternative argument for proving that solutions
$\bigtriangleup(d_i)$ exist and are at most two, finiteness being
rather obvious.
With $d_i\in v(P_1)\approx P_1$ as unknowns, the determinantal
condition on $(d_i,d_j)\in (P_1)^2$ expressing collinearity with $p_{ij}$
is of type $(2,2)$ but contains the diagonal as improper solutions. Thus
we actually have a $(1,1)$ condition. On $(P_1)^3$ we intersect accordingly
three equations: of type $(1,1,0),\ (0,1,1)$ and $(1,0,1)$. This yields,
counting multiplicity, two solutions.

\medskip \noindent
For the more precise statement in our lemma, we observe that, in the case of a proper secant, there's an
involution of $P_2$, induced from an involution of $P_1\approx v(P_1)$, and  which keeps $\ell$ pointwise fixed.
(One extends to $P_2$ the involution of $P_1=v(P_1)$ fixing the intersection with $\ell$.) Thus, a triangle
solution produces a `reflected' second solution.

\medskip \noindent
When $\ell$ is tangent, this is no longer the case. Indeed, if we would have two solutions:
$\bigtriangleup(d_i)$ and $\bigtriangleup(d'_i)$, there would be an involution of $P_1\approx v(P_1)$ taking one
onto the other, defined by tracing lines through the perspective point $p$ and exchanging the two intersection
points with the conic. The associated transformation of $P_2=P(Sym_C(2))$ would have to fix $\ell$ pointwise,
since it must fix $p_{ij}$, but it also fixes the line through the tangency points of the two tangents from $p$
to the conic. This is a contradiction and the lemma is proven. \ \ \ $\Box$

\medskip \noindent
To conclude this analysis, we do some dimension counts.

\medskip \noindent
The quasi-projective variety:

$$ B^3_{21}=\{ (q_1,q_2,q_3,b)\ :\ b\ \mbox{is a common basket for} \ q_1, q_2 \
\mbox{and}\  q_3 \  \}$$

\noindent projects to $P_9$ (the closure of the space of possible baskets $b$) by $(q_1,q_2,q_3,b)\mapsto b$.
The fibers are irredicible open subsets of the third Cartesian power of the cone from $b$ over the rank one
locus. This gives irreducibility and dimension:\ $9+3\times 4=21$ for $B^3_{21}$.

\medskip \noindent
From the perspective of our characterization, $B^3_{21}$ is obtained
(up to birational equivalence) by running a plane through a generic pencil
of cones with fixed vertex, and then choosing a triangle in this
plane with edges passing through the three rank two points. This gives
dimension 21 as $11+7+3$.

\medskip \noindent
Similar counts can be performed for various subvarieties
of $B^3_{21}$. For example, when our pencil degenerates to one
contained in ${\cal R}^2_6$ (i.e. becomes a point of ${\cal F}^2_6$),
we are generically in the case addressed by our previous lemma, and
the dimension of the corresponding subvariety is:\ $6+7+3\times 2=19$.

\medskip \noindent
We can regroup now the main results in this section in the
form of a complement to Proposition 0.3 :

\begin{prop}
Consider the quasi-projective variety

$$ B^3_{21}=\{ (q_1,q_2,q_3,b)\ :\ b\ \mbox{is a common basket for} \
q_1, q_2 \ \mbox{and}\  q_3 \  \} \subset (P_9-{\cal R}^2_6)^4 $$

\medskip \noindent
One can distinguish several closed subvarieties, according to the diagram:

$$
\begin{array}{ccccccc}
B^3_{21} & \supset & C^3_{19} & \supset & E^3_{15} & \supset &
F^3_{14} \\
\cup & \  & \cup  & \  & \cup & \  & \  \\
H^3_{14} & \subset  & D^3_{18}  & \  & G^3_{13} & \ & \
\end{array}
$$

\medskip \noindent
These subvarieties are determined by the following geometrical conditions
on the quadrics $q_i,\ i=1,2,3$:

\medskip \noindent
$(B^3_{21}):$\ the triple $(q_1,q_2,q_3)$ of a generic point
$(q_1,q_2,q_3,b)\in B^3_{21}$ spans a plane $[q_1,q_2,q_3]=P_2$
containing a pencil of singular quadrics with exactly three
distinct rank two points $p_i\in [q_j,q_k]$ (and rank three elsewhere);

\medskip \noindent
$(C^3_{19}):$ \ the triple $(q_1,q_2,q_3)$ of a generic point
$(q_1,q_2,q_3,b)\in C^3_{19}$ spans
a plane $[q_1,q_2,q_3]$ containing a secant or tangent of the rank one
locus (i.e. a line from the family ${\cal F}^3_6$);

\medskip \noindent
$(D^3_{18}):$ \ the triple $(q_1,q_2,q_3)$ of any point
$(q_1,q_2,q_3,b)\in D^3_{18}$ spans
a plane $[q_1,q_2,q_3]$ containing a proper secant of the rank one
locus, and one edge $[q_i,q_j]$ passes through a rank one point
on this secant;

\medskip \noindent
$(E^3_{15}):$ \ the triple $(q_1,q_2,q_3)$ of  a generic point
$(q_1,q_2,q_3,b)\in E^3_{15}$ spans
a proper secant plane of the rank one locus;

\medskip \noindent
$(F^3_{14}):$ \ the triple $(q_1,q_2,q_3)$ of  a generic point
$(q_1,q_2,q_3,b)\in F^3_{14}$ spans
a proper secant plane of the rank one locus and the triangle
$\bigtriangleup(q_i)$ is in perspective with the triangle of
double-planes in that span;

\medskip \noindent
$(G^3_{13}):$\ the triple $(q_1,q_2,q_3)$ of any point
$(q_1,q_2,q_3,b)\in G^3_{13}$ spans a pencil
of conics with three rank two points, and condition $(c1)$ is
satisfied on the pencil for some ordering of these points.

\medskip \noindent
$(H^3_{14}):$ \ the triple $(q_1,q_2,q_3)$ of any point
$(q_1,q_2,q_3,b)\in H^3_{14}$ spans a pencil
which meets the rank one locus.

\medskip \noindent
Any point of $B^3_{21}$ satisfies one of the conditions above i.e. it is
either already generic on $B^3_{21}$, in the specified sense, or is to be
found on (at least) one of these closed subvarieties.

\medskip \noindent
The dimension of the subvariety is indicated by the subscript, and
the number of points $(q_1,q_2,q_3,b)$ with the same projection
$(q_1,q_2,q_3)$ {\em for some open dense set } in
each subvariety is tabulated below:

\begin{center}
\begin{tabular}{ccccccc}
$B^3_{21}$ & $C^3_{19}$ & $D^3_{18}$ & $E^3_{15}$ & $F^3_{14}$  &
$G^3_{13}$ & $H^3_{14}$ \\ \hline
1        & 2       & 1              & 6        & 7           & 1  & $P_1$
\end{tabular}
\end{center}

\end{prop}

\medskip \noindent
{\bf Remark:} \ All triples involved in the families
$E^3_{15},F^3_{14},G^3_{13}$
are made of cones with common vertex, and relate to issues
of tangency for conics. They are less relevant  for questions
of common tangents to quadric surfaces because all lines through
the common vertex are always common tangents.

\section{Four quadrics in a basket}

\medskip \noindent
In this section we prove and expand Proposition 0.4.

\medskip \noindent
Let $q_i,\ i=1,...,4$ be four quadrics with a common basket $b$, and let $d_i$ be the double-plane in the pencil
$[b,q_i]$. {\em Generically}, the $q_i$'s would span a three-space $P_3\approx [q_1,...,q_4]$, the $d_i$'s would
span another three-space $[d_1,...,d_4]$, and the two would meet in a plane $P_2\approx [q_1,...,q_4]\cap
[d_1,...,d_4]$. This plane contains a {\em complete quadrilateral} made of the four lines
$\ell_i=[q_j,q_k,q_l]\cap [d_j,d_k,d_l]$. These lines are pencils of cones and belong to the family ${\cal
F}^3_{11}$. Clearly the pencil $[q_i,q_j]$ passes through the intersection point $p_{kl}=\ell_k\cap
\ell_l=[q_i,q_j]\cap [d_i,d_j]$, and the six points $p_{kl}$ are all rank two points.

\medskip \noindent
We begin our analysis from this end, and establish some facts about
complete quadrilaterals.

\medskip \noindent
{\bf Definition:}\ A {\em complete quadrilateral} consists of
four lines in general position in $P_2$, that is: no three are
concurrent.

\medskip \noindent
Equivalently, a {\em complete quadrilateral} is a projective planar configuration $(6_2,4_3)$ of six points and
four lines, with every point incident to two lines, and every line incident to three points. The above
labelling: \ $p_{ij}=\ell_i\cap \ell_j, \ \{i,j\}\subset \{1,2,3,4\}$ will be normally adopted.

\begin{lemma}
{\bf (Cayley)}\  Every complete quadrilateral can be obtained by
intersecting the faces of a tetrahedron in $P_3$ with a plane $P_2$
avoiding its vertices.
\end{lemma}

\medskip \noindent
{\bf Remark:}\ It follows from Desargues' theorem that if two
tetrahedra in $P_3$ cut in this fashion the same complete quadrilateral
in a $P_2$, then they are in perspective, that is: with proper labelling,
the four lines through corresponding vertices meet in the same point
(the perspective point). \ \ \ $\Box$

\begin{lemma}
Suppose we have a complete quadrilateral in $P_2\subset P_9$, made of pencils $\ell_i\subset {\cal R}^3_8$, with
rank two points $p_{ij}=\ell_i\cap \ell_j$ and rank three  elsewhere.

\medskip \noindent
Then each $\ell_i$ is a pencil of cones with fixed vertex $v_i$, and either:

\medskip
i)\ \  the $v_i$'s are in general position, or:

ii)\ all $v_i$'s coincide.

\medskip \noindent
In the former case, if we denote by $d_i$ the double-plane supported
by the span $P_2\approx [v_j,v_k,v_l]\subset P_3$, we obtain the
{\em unique tetrahedron of rank one points} which contains the initial
$P_2\subset P_9$ in its span, and hence produces the given complete
quadrilateral by the four traces of its faces.

\medskip \noindent
The latter case, when $P_2\subset P_5=P(Sym_C(3))\subset P_9=P(Sym_C(4))$,
is addressed in the next lemma.
\end{lemma}

\medskip \noindent
{\em Proof:} \ By Lemma 3.1, each $\ell_i$ is a pencil of cones with fixed vertex $v_i\in P_3$. We look at the
traces of these pencils on a plane $P_2\subset P_3$ chosen away from the vertices.

\medskip \noindent
 $\ell_i$
traces a pencil of conics passing through {\em four points in general
positions}, and the three rank two points on it $p_{ij}, \ j\neq i$
correspond to the three pairs of lines, with complementary pairs
of points on it.

\medskip \noindent
These three pairs of rank two conics have consequently {\em non-collinear} singularities $s_{ij}$. The line
$[v_i,s_{ij}]$ is clearly the singularity axis for $p_{ij}$ and contains therefore $v_j$. Thus, two of the
vertices $v_j,\ v_k$ cannot coincide without being both equal to $v_i$, and repeating this argument shows that
we can either have:

\medskip
i)\ \  $v_i,\ i=1,...,4$ in general position in $P_3$, or:

ii)\ all $v_i$'s equal.

\medskip \noindent
We pursue here the first case, and adopt as projective coordinates
$(x_1:...:x_4)\in P_3$ those corresponding to the reference tetrahedron
$v_i$, that is: $\{ x_i=0\}$ is the face $[v_j,v_k,v_l]$.

\medskip \noindent
Then, the quadrics in the pencil $\ell_i$ do not involve the variable $x_i$, and so $p_{ij}$ involves neither
$x_i$, nor $x_j$.

\medskip \noindent
But $\ell_i=\lambda p_{ij}+\mu p_{ik}$ contains $p_{il}$, which has neither $x_i$, nor $x_l$, and this can
happen only when $p_{ij}$ has no $x_kx_l$ term, and $p_{ik}$ has no $x_jx_l$ term. This means: \ $\ell_i\subset
[x_j^2,x_k^2,x_l^2]$, hence our complete quadrilateral is in the span of $d_i=x_i^2,\ i=1,...,4$.

\medskip \noindent
Uniqueness follows from the uniqueness part in Lemma 3.1.\ \ \ $\Box$

\medskip \noindent
{\bf Remark:}\ Degree considerations imply that a complete quadrilateral made of singular quadrics cannot have a
line contained in ${\cal R}^2_6$ unless its plane is contained in ${\cal R}^3_8$.

\medskip \noindent
When all vertices coincide, projection from the common vertex
reduces the problem to the space of conics:

\begin{lemma}
Suppose we have a complete quadrilateral in $P_2\subset P_5=P(Sym_C(3))$, made of pencils of conics $\ell_i$,
with rank two points $p_{ij}=\ell_i\cap \ell_j$ and rank three somewhere. Suppose further that $P_2$ is not a
secant plane of the rank one locus i.e. not the span of three rank one points.

\medskip \noindent
Then, up to relabelling, there's a unique tetrahedron of rank one points
$d_i,\ i=1,...,4$,
which contains the initial plane in its span, and produces the
complete quadrilateral by the four traces of its faces.
\end{lemma}

\medskip \noindent
{\em Proof:}\ The rank stratification in the space of conics reads:

$$ P_5=P(Sym_C(3))\supset {\cal R}^2_4\supset {\cal R}^1_2=v(P_2) $$

\noindent
with ${\cal R}^2_4$ of degree three, and ${\cal R}^1_2$ of degree four.

\medskip \noindent
Thus, under our assumptions, $P_2\cap {\cal R}^2_4$ is a {\em cubic curve}
(but not degenerated into three lines),
with six distinct points $p_{ij}$ on it, subject to four collinearity
conditions.

\medskip \noindent
The assumption of some rank three point means we can apply the argument in Lemma 3.1 to one of the lines, say
$\ell_4$, and find a unique trio of rank one conics $d_1,d_2,d_3$ whose span contains $\ell_4$.

\medskip \noindent
Because of our other assumption,
$P_2$ and $[d_1,d_2,d_3]$ span a three-space $P_3$,
which must meet the rank one locus in at least one more point,
distinct from the other three. However, there might be a whole
conic of rank one points in the intersection.

\medskip \noindent
The case of a single new rank one point $d_4$ obviously corresponds
to a complete quadrilateral with no line in ${\cal R}^2_4$, and is
thus resolved by the quadruple $d_i,\ i=1,...,4$.

\medskip \noindent
The alternative requires one line of the complete quadrilateral to be contained in ${\cal R}^2_4$, and if we
call it $\ell_1$, then, with the natural labelling in use, the Veronese conic $v(P_1)$ introduced in the
intersection will pass through $d_2$ and $d_3$. Thus, $\ell_1$ is a secant or tangent to this conic.

\medskip \noindent
The restriction of ${\cal R}^2_4$ to our $P_3$ consists therefore of the plane $P_2\supset v(P_1)$ and the {\em
cone} from $d_1$ over $v(P_1)$. (Indeed, $d_1$ must be singular on the residual quadric.)

\medskip \noindent
Thus, in order to find our $d_4$ we simply intersect $[p_{23},d_1]$
with $v(P_1)$.

\medskip \noindent
This yields the desired quadruple of
rank one points, clearly unique up to relabelling. \ \ \ $\Box$

\medskip \noindent
We need to investigate also the most degenerate case, when the entire
plane of the complete quadrilateral lies in ${\cal R}^2_4$.
Since our lines must be pencils with fixed singularity, we are actually
envisaging the case $P_2=P(Sym_C(2))\supset {\cal R}^1_1= v(P_1)$.

\medskip \noindent
Thus, the question is: \ given a complete quadrilateral in $P_2=P(Sym_C(2))$,
when is there a quadruple of rank one points $d_i\in v(P_1)$ such that
$p_{ij}\in [d_k,d_l],\ \{i,j\}\cup \{k,l\}=\{1,2,3,4\} $ ?

\medskip \noindent
The space of complete quadrilaterals in $P_2$ is birationally equivalent to $(P_2)^4$, hence eight-dimensional.
On the other hand, the space of ordered quadruples of double-points is $(v(P_1))^4\approx (P_1)^4$, and a given
quadruple has the required relation with a two-parameter family of complete quadrilaterals. Thus, only a six
dimensional subfamily of complete quadrilaterals can be `solved' in this sense, and we need a `codimension two'
condition $(c2)$ satisfied.

\medskip \noindent
In order to streamline some of our statements, we introduce:

\medskip \noindent
{\bf Definition:}\ {\em (Typical, special, and solvable complete
quadrilaterals)}

\medskip \noindent
A {\em complete quadrilateral} in $P_9=P(Sym_C(4))$
will be called {\bf typical} when defined by the traces of the
four faces of a tetrahedron with vertices at rank one points on a
sectioning plane avoiding these vertices.

\medskip \noindent
A {\em complete quadrilateral} will be called {\bf special} when
contained in the span of a proper secant plane of the rank one
locus (i.e. a plane with exactly three rank one points), and has
one vertex of rank one, with the remaining five of rank two.

\medskip \noindent
A {\em complete quadrilateral} will be called {\bf solvable} when
contained in the span of a conic of rank one points (i.e. a plane
$P_2\supset v(P_1)$), and  there's a quadruple of rank one points
$d_i\in v(P_1)$, such that the six vertices satisfy:\
$p_{ij}\in [d_k,d_l]$. \ \ \ $\Box$

\begin{lemma}
Given a {\em solvable} complete quadrilateral in $P_2=P(Sym_C(2))$,
there's a {\em unique} corresponding solution $d_i,\ i=1,...,4$.
\end{lemma}

\medskip \noindent
{\em Proof:}\ Two solutions $d_i$ and $d'_i$ would be necessarily in perspective. As in Lemma 3.2, intersecting
the conic with lines through the perspective point and exchanging the two intersection points defines an
involution of $P_1\approx v(P_1)$ which exchanges the two solutions. However, the associated involution of $P_2$
would have to fix the complete quadrilateral and thus be the identity. The contradiction proves the claim. \ \ \
$\Box$

\begin{lemma}
The plane of a {\em typical complete quadrilateral} belongs
to one of the following (disjoint) families:

\medskip \noindent
$(\Phi_{15}):$\ planes $P_2\subset P_9$ where ${\cal R}^3_8=\{ det(Q)=0\}$
restricts to four lines in general position (the complete quadrilateral
itself);

\medskip \noindent
$(\Psi_{12}):$\ planes $P_2\subset P_5\subset {\cal R}^3_8$, with $P_5=P(Sym_C(3))$ in adequate coordinates,
where (relative to the corresponding space of conics) ${\cal R}^2_4$ restricts to a cubic curve other than three
lines.

\medskip \noindent
Every complete quadrilateral supported in a plane of type $(\Psi)$,
and with vertices at rank two points, is admissible.

\medskip \noindent
The tetrahedron of rank one points defining an admissible complete
quadrilateral is {\em unique}.
\end{lemma}

\medskip \noindent
{\em Proof:}\ In view of the above discussion, all that remains to
be shown is that for type $(\Phi)$, the lines of singular quadrics
meet at rank two points. But this follows from our investigation of
lines in ${\cal R}^3_8$ of section 1.

\medskip \noindent
Indeed, the lines in $P_2\cap {\cal R}^3_8$ must belong to one of
the families: \ ${\cal F}^3_{11}-{\cal F}^2_6$ or
${\cal M}^3_{11}-{\cal M}^2_7$. However, in the latter case, the
pencil would have a single rank two point. Yet, there are three
points on the pencil which are singularities of the restricted
determinantal quartic, namely: its intersections with the other three
lines. Since our pencil cannot be tangent to ${\cal R}^3_8$ in
more than two points, this latter case must be discarded.

\medskip \noindent
This leaves us with a pencil of type $(F)$ i.e. a pencil of cones with fixed vertex, and the three rank two
points on such pencils must be, in our case, the intersections with the other three lines. \ \ \ $\Box$

\medskip \noindent
{\bf Note:}\ {\em While the definition of a {\em typical, special, or
solvable complete quadrilateral} uses a quadruple of rank one points,
we have seen above that these properties can be detected directly  from
the complete quadrilateral itself, and depend essentially on the
position of its span with respect to the rank stratification of
$P_9=P(Sym_C(4))$. The quadruple of rank one points can be `reconstructed'
from this type of information. Even without searching here for the
explicit form of condition $(c2)$, we shall use henceforth: {\em
typical, special, or solvable} in the sense of {\em a property which
needs no explicit mention of four rank one points}.}

\medskip \noindent
We can list now the possible types of configurations for $(q_i)_i,
b,(d_i)_i, i=1,2,3,4$\ :\ four quadrics, a common basket, and rank one
points $d_i\in [b,q_i]$. We shall see, in the spirit of the above
note, that these classes need no explicit mention of $b$ and $(d_i)_i$.

\medskip \noindent
First, for a {\bf three-dimensional span}:\ $[q_1,q_2,q_3,q_4]=P_3$.

\medskip \noindent
(B):\ $[d_1,...,d_4]=P_3$

\medskip \noindent
(C):\ $[d_1,...,d_4]=P_2$ because $b\in [q_i,q_j]$ and hence $d_i=d_j$

\medskip \noindent
(D):\ $[d_1,...,d_4]=P_2$ because $d_i\in v(P_1),\ i=1,...,4$

\medskip \noindent
Next, for a {\bf two-dimensional span}:\ $[q_1,q_2,q_3,q_4]=P_2$.

\medskip \noindent
(E):\ $[b,q_1,...,q_4]=[d_1,...,d_4]=P_3$

\medskip \noindent
(F):\ $b\in [q_i,q_j],\ d_{ij}=d_i=d_j,\  [d_{ij},d_k,d_l]=P_2$

\medskip \noindent
(G):\ $[b,q_1,...,q_4]=P_3,\  [d_1,...,d_4]=P_2$

\medskip \noindent
(H):\ $b\in [q_i,q_j]\cap [q_k,q_l]$

\medskip \noindent
(I):\ $b\in [q_i,q_j,q_k]=P_1$

\medskip \noindent
Lastly, for a {\bf one-dimensional span}:\ $[q_1,q_2,q_3,q_4]=P_1$.

\medskip \noindent
(J):\ $b\in [q_1,...,q_4]=P_1$

\medskip \noindent
This list structures our extension of Proposition 0.4. In order to
indicate inclusion, rather than adjacency, we consider our subvarieties
as {\em closed subvarieties of} $B^4_{25}$, that is: as the closure in
$B^4_{25}$ of the locus described by some generic property.

\begin{prop}
Consider the quasi-projective variety:

$$ B^4_{25}=\{ (q_1,q_2,q_3,q_4,b) \ : \ b\ \mbox{is a common basket for}\
q_i, \ i=1,...,4 \  \} \subset (P_9-{\cal R}^2_6)^5 $$

\medskip \noindent
One can distinguish closed subvarieties:

$$
\begin{array}{ccccccc}
F^4_{17} & \subset & C^4_{22} & \       & H^4_{19} & \  & \  \\
\cap     & \       & \cap     & \       & \cap     & \  & \  \\
E^4_{18} & \subset & B^4_{25} & \supset & G^4_{20} & \supset &
J^4_{15}=H^4_{19}\cap I^4_{19} \\
\        & \       & \cup     & \       & \cup     & \  & \  \\
\        & \       & D^4_{21} & \       & I^4_{19} & \  & \
\end{array}
$$

\noindent
according to the following geometrical conditions:

\medskip \noindent
$(B^4_{25}):$ \ the quadrics $(q_i)$ of a generic point
$(q_1,...,q_4,b)\in B^4_{25}$ define
a reference tetrahedron in a three-space, tracing a {\em typical
complete quadrilateral} on some two-subspace avoiding its vertices;

\medskip \noindent
$(C^4_{23}):$\ the quadrics $(q_i)$ of a generic point
$(q_1,...,q_4)\in C^4_{23}$ define
a reference tetrahedron in a three-space, tracing a {\em special
complete quadrilateral} on some two-subspace avoiding its vertices;

\medskip \noindent
$(D^4_{21}):$\ the quadrics $(q_i)$ of any point
$(q_1,...,q_4)\in D^4_{21}$ span a
three-space containing a $P_2\supset v(P_1)$, and the faces $[q_i,q_j,q_k]$
trace a {\em solvable complete quadrilateral} in this $P_2$;

\medskip \noindent
$(E^4_{18}):$\ the quadrics $(q_i)$ of a generic point
$(q_1,...,q_4,b)\in E^4_{18}$ span
a plane containing a {\em typical complete quadrilateral} with vertices
$p_{ij}\in [q_k,q_l]$;

\medskip \noindent
$(F^4_{17}):$\ the quadrics $(q_i)$ of a generic point
$(q_1,...,q_4,b)\in F^4_{17}$ spans
a proper secant plane of the rank one locus and contains a {\em special
complete quadrilateral} with vertices $p_{ij}\in [q_k,q_l]$;

\medskip \noindent
$(G^4_{20}):$\ the quadrics $(q_i)$ of a generic point $(q_1,...,q_4,b)\in G^4_{20}$ span a plane containing a
secant or a tangent $\ell$ of the rank one locus, and the four quadrics $q_i$ lie on a conic which contains the
two rank one points $A$ and $B$ of the secant, or is tangent to $\ell$ at $T=A=B$ when $\ell$ is a tangent;

\medskip \noindent
$(H^4_{19}):$\ the quadrics $(q_i)$ of a generic point
$(q_1,...,q_4,b)\in H^4_{19}$ span
a plane, and for some $\{i,j,k,l\}=\{1,2,3,4\}$, both $[q_i,q_j]$ and
$[q_k,q_l]$ meet the rank one locus;

\medskip \noindent
$(I^4_{19}):$\ the quadrics $(q_i)$ of a generic point
$(q_1,...,q_4,b)\in I^4_{19}$ span a plane; three of them are
on a line meeting the rank one locus, and the span contains
another double-plane;

\medskip \noindent
$(J^4_{15}):$\ the quadrics $(q_i)$ of any point
$(q_1,...,q_4,b)\in J^4_{15}$ span a line
meeting the rank one locus.

\medskip \noindent
Any quadruple of quadrics allowing a common basket enters one of
the configurations listed above as $(B)$ to $(J)$, and satisfies
the corresponding property stated above.
\end{prop}

\medskip \noindent
{\em Proof:}\ Since almost all relevant arguments have already been
presented, we fill in a few remaining details.

\medskip \noindent
For the $G^4_{20}$ family, the conic through the four points $q_i$ which also contains the rank one points $A$
and $B$ of the secant (or is tangent at $T=A=B$ to $\ell$ in case of a tangent), is obviously the projection
from $b$ to $[q_1,...,q_4]=P_2$ of the Veronese conic $v(P_1) \subset P_2 \supset \ell$.

\medskip \noindent
The existence of this conic is equivalent to conditions we label $(c3)$ for the six points $p_{ij}$ on the
secant or tangent. These conditions reflect the fact that $p_{ij}=\ell\cap [d_k,d_l]$ and $d_i$ are on the conic
$v(P_1)$.

\medskip \noindent
In case $\ell$ is a secant, $A$ and $B$ will be the two points of $\ell\cap v(P_1)$. Considering the projections
of the conic $v(P_1)$ from $d_i$, respectively $d_j$, onto $\ell$, we obtain:

$$ (A,B;p_{jk},p_{jl})=(A,B;p_{ik},p_{il}), \ \ \  \{i,j,k,l\}=\{1,2,3,4\} $$

\noindent
and we call $(c3)$ the collection of these cross-ratio relations.

\medskip \noindent
In case $\ell$ is tangent, we have $A=B=T=\ell\cap v(P_1)$. Again, projecting from $d_i$, then $d_j$, we obtain:

$$ (T,d_j;d_k,d_l)=(T,p_{kl};p_{lj},p_{jk}) $$

$$(d_i,T;d_k,d_l)=(p_{kl},T;p_{li},p_{ik}) $$

\noindent
with the first cross-ratio on the conic, and the second on the tangent.
We can eliminate the $d_i$'s from the resulting system by first taking
the product of the two left-hand sides above with $(d_i,d_j;d_k,d_l)$,
to get:

$$ 1=(T,p_{kl};p_{lj},p_{jk})\cdot (p_{kl},T;p_{li},p_{ik})\cdot
(d_i,d_j;d_k,d_l) $$

\noindent
and hence:

$$(T,p_{kl};p_{lj},p_{jk})\cdot (p_{kl},T;p_{li},p_{ik})=
(T,p_{ij};p_{jl},p_{li})\cdot (p_{ij},T;p_{jk},p_{ki}) $$

\medskip \noindent
The condition $(c3)$ will be the collection of these equations, should we
have tangency, and not a proper secant.

\medskip \noindent
One verifies that, for $p_{ij}$ of rank two, $(c3)$ is sufficient for
finding $d_i$'s on the conic
with $p_{ij}\in [d_k,d_l]$. By the same argument as in Lemma 3.2, there
are two solutions in the secant case, and one solution in the tangent
case. \ \ \ $\Box$

\section{A double-four example}

\medskip \noindent
In this section we study a particular configuration, made of two
groups of four quadrics. Each group has {\em linearly independent and
smooth quadrics}, and the four in one group are common baskets for
the four in the other, hence the designation ``double-four''.

\medskip \noindent
It turns out that, up to projective transformations of $P_3$, this configuration is {\em uniquely determined} by
the stated property. Our proof uses the criteria developed above, in particular the configuration of points and
lines created by four mashed complete quadrilaterals. It will be identified as the {\em Reye configuration}
$(12_4,16_3)$. \cite{HC-V} \cite{Dol}

\medskip \noindent
We start with the eight `diagonal' quadrics:\ $\pm x_1^2\pm x_2^2\pm
x_3^2\pm x_4^2$. The two groups of four quadrics are those of
positive, respectively negative determinant.

\medskip \noindent
Obviously permutations and changes of two signs preserve the two
groups, while changing one sign exchanges the groups. Thus,
our labelling here is mostly a matter of convenience. We put:

$$ q_i=(\sum_{j=1}^4 x_j^2)-2x_i^2,\ \ d_i=x_i^2, \ \ \ i=1,...,4 $$

$$ b^1=\sum_{j=1}^4 x_j^2, \ b^{\beta}=(\sum_{j=1}^4 x_j^2)-2(x_{\beta}^2+
x_{\beta-1}^2),\ \ \ \beta=2,3,4 $$

\noindent Clearly, we have three tetrahedra, spanning the same three-space:

$$ P_3=[q_1,...,q_4]=[b^1,...,b^4]=[d_1,...,d_4] $$

\medskip \noindent
Intersecting the faces of the first two tetrahedra yields sixteen
lines:

$$ \ell_i^{\alpha}=[q_j,q_k,q_l]\cap [b^{\beta},b^{\gamma},b^{\delta}] $$

\noindent
and the {\em same} collection of sixteen lines obtains from intersecting
the faces of the first and last tetrahedra, or second and last. Thus,
there are twelve planes, each containing four lines, with each line
contained in three planes.

\medskip \noindent
This is obviously the {\em dual} of a $(12_4,16_3)$ configuration, but
we would rather distinguish a direct $(12_4,16_3)$ configuration by
taking into account the {\em twelve rank two points} which lie on
the sixteen lines, with three points on each line, and each point
incident to four lines.

\medskip \noindent
Our labelling is now going to show its bias for the $q_i$ and $d_j$
tetrads, but one should remain aware of the perfectly equivalent
role of the third tetrad $b^{\alpha}$. We put:

$$ p_{ij}^+=x_k^2+x_l^2=[d_k,d_l]\cap [q_i,q_j] $$

$$ p_{ij}^-=x_k^2-x_l^2=[d_k,d_l]\cap [q_k,q_l] $$

\noindent
Thus, each edge of our three tetrahedra has exactly two of the twelve
rank two points $p_{ij}^{\pm}$.

\medskip \noindent
One can identify now the collection of points and lines $(p_{ij}^{\pm}, \ell_k^{\alpha})$ with their incidence
relations, as the {\em Reye configuration} $(12_4,16_3)$. The latter is usually depicted in some affine (real)
part $R^3\subset P_3$ by means of a cube, with $p_{ij}^{\pm}$ corresponding with the eight vertices of the cube,
its center, and the three points at infinity determined by the three distinct directions of the edges. The
sixteen lines $\ell_k^{\alpha}$ are the twelve edges and the four diagonals of the cube. \cite{HC-V}\cite{Dol}

\medskip \noindent
Our three tetrahedra $(q_i)_i, (d_j)_j, (b^{\alpha})_{\alpha}$, should
be seen in this model as the three tetrads of planes determined, each,
by two opposite faces of the cube together with the two diagonal
planes of the cube which cut the two diagonals in these faces.

\medskip \noindent
We are going to see this configuration emerging from any ``double-four'',
and obtain:

\begin{theorem}
Suppose $(q_i)_i$ and $(b^{\alpha})_{\alpha}$ form a double-four
configuration of {\em smooth quadrics} in $P_3$, that is:

$q_i$ are linearly independent and are common baskets for
$(b^{\alpha})_{\alpha}$, and

$b^{\alpha}$ are linearly independent and common baskets for $(q_i)_i$.

\medskip \noindent
Then, there is a projective automorphism of $P_3$ which carries
$(q_i)_i$ and $(b^{\alpha})_{\alpha}$ to the standard double-four
configuration presented above.
\end{theorem}

\medskip \noindent
The proof requires a number of lemmas. Throughout, we let $q_i,b^{\alpha}$
stand for a double-four as in the theorem, but the notation should be
understood as separate from the one used in describing the standard
double-four.

\medskip \noindent
Since $b^{\alpha}$ is a basket for $d_i$, the pencil $[b^{\alpha},q_i]$
has a (unique) rank one point $d_i^{\alpha}$. Considering $(q_i)_i$
as a quadruple with common basket $b^{\alpha}$, linear independence
and Proposition 4.6. produce an associated {\em complete quadrilateral}:

$$ \ell_i^{\alpha}=[q_j,q_k,q_l]\cap
[d_j^{\alpha},d_k^{\alpha},d_l^{\alpha}], \ \ \ i=1,...,4 $$

\noindent
with six vertices $p_{ij}^{\alpha}=\in [q_k,q_l]$ at points of rank
at most two.

\medskip \noindent
However, an edge $[q_i,q_j]$ has at most two points of rank two, or
a single rank one point and no rank two. Thus, the collection
$p_{ij}^{\alpha},\ \{i,j\}\subset \{1,2,3,4\}, \alpha=1,...,4$
consists of at most twelve {\em distinct} points. We'll see that
there must be precisely twelve distinct points, all of rank two.

\medskip \noindent
The fact that there can be no rank one point $p_{ij}^{\alpha}$
follows from the observation that if one complete quadrilateral
were special, all would be special, with the same rank one point
on some edge $[q_i,q_j]$. But this would force all baskets on
that same edge, contradicting linear independence.

\begin{lemma}
The sixteen lines $\ell_i^{\alpha}$ are all distinct and none lies in ${\cal R}^2_6$.
\end{lemma}

\medskip \noindent
{\em Proof:}\
Since the four complete quadrilaterals corresponding to $b^{\alpha},\
\alpha=1,...,4$ must be distinct, their planes $\pi_{\alpha}$ are
distinct.

\medskip \noindent
Thus, two complete quadrilaterals can share at most one line, and if they do share one line, it must be in
${\cal R}^2_6$, for otherwise the three rank two points on it would be the same for the two quadrilaterals, and
this would already force the two baskets to coincide.

\medskip \noindent
We can make now a first estimate of how many of the points
$p_{ij}^{\alpha}$ should be {\em distinct}. A first quadrilateral brings
in six, a second at least three more, and a third at least one
more. We find a {\em minimum} of ten. Thus, at least four of the
edges $[q_i,q_j]$ must meet the rank two locus in two distinct points.

\medskip \noindent
This already rules out the possibility of a solvable complete quadrilateral
amongst our tetrad. Indeed, we would then have a Veronese conic
$v(P_1)\subset P_2\subset [q_1,...,q_4]=P_3$, and the restriction of
${\cal R}^3_8$ to this $P_3$  would decompose into the plane
$P_2\subset {\cal R}^2_6$ of the conic, counted twice, and a quadric.

\medskip \noindent
This residual quadric must carry more than one rank two point away from the double $P_2$, otherwise four edges
$[q_i,q_j]$ would be concurrent. But two singularities make the quadric a two-plane. And this is not good enough
to allow four edges of the $q_i$ tetrahedron to intersect the rank two locus in two points, unless the two-plane
is actually a double-plane made of points of rank at most two. But this forces $[q_1,...,q_4]$ entirely into
${\cal R}^3_8$: a contradiction.

\medskip \noindent
Thus, {\em all four complete quadrilaterals are typical}.

\medskip \noindent
Let us return to the hypothesis that two of them share a line, say
$\ell_i^{\alpha}=\ell_i^{\beta}=\pi_{\alpha}\cap \pi_{\beta}\subset {\cal R}^2_6$. The restriction of ${\cal
R}^3_8$ to $P_2=[q_j,q_k,q_l]$ then decomposes into the double-line $\ell_i^{\alpha}$ and a conic. Since at
least one edge in our face requires a second rank two point on it, the conic must be a two-line. We would arrive
at a contradiction, as in the previous argument, should we know that two edges require two rank two points.

\medskip \noindent
But now we can review the estimate of a minimum of ten distinct points
$p_{ij}^{\alpha}$, and see it based on repeated common lines. Hence,
either there are at least eleven distinct points $p_{ij}^{\alpha}$ and
we get our contradiction, or we find two lines contained both in
${\cal R}^2_6$ and some face $[q_i,q_j,q_l]$, a contradiction again.

\medskip \noindent
This proves that our four complete quadrilaterals cannot have
lines in common. With that, the estimate on the cardinality of
the $p_{ij}^{\alpha}$ collection is lifted to the maximum twelve:
six from a first quadrilateral, at least four more from a second,
and at least two more from a third. Moreover, any two complete
quadrilaterals share {\em exactly two such points} and no three
can share the same point.

\medskip \noindent
With the established fact that all six edges $[q_i,q_j]$ carry two rank two points, we conclude as above that no
line $\ell_i^{\alpha}$ is contained in ${\cal R}^2_6$. \ \ \ $\Box$

\medskip \noindent
The proof has given at the same time:

\begin{cor}
There are exactly twelve distinct points in the set:

$$ \{ p_{ij}^{\alpha} \ : \ \{ i,j\}\subset \{ 1,2,3,4\}, \ \alpha=1,..,4\} $$

\noindent
with two of them on each edge $[q_i,q_j]$.

\medskip \noindent
There are three of them on each of the sixteen lines $\ell_i^{\alpha}$, and four of these lines are passing
through each point. Thus, $( p_{ij}^{\alpha},\ell_k^{\beta})$ defines a $(12_4,16_3)$ configuration.
\end{cor}

\medskip \noindent
Our aim now is to prove that the span of the $q_i$'s is the same as the
span of the $b^{\alpha}$'s, and the same as the span of the
$d_i^{\alpha}$'s, which turn out in fact to be just four distinct
rank one points with symmetric role towards the two other tetrads.

\medskip \noindent
In order to simplify statements, we shall refer to the {\em twelve
points of rank two} $p_{ij}^{\alpha}$ as {\em marked points}.

\begin{lemma}
For any edge $[q_i,q_j]$ and marked rank two point on it, there's some
edge $[b^{\alpha},b^{\beta}]$ meeting it at that point.
\end{lemma}

\medskip \noindent
{\em Proof:}\ Each basket $b^{\alpha}$ is associated with one
of the two marked points on $[q_i,q_j]$, hence we'll find
two baskets, say $b^{\alpha}$ and $b^{\beta}$, with
$p_{kl}^{\alpha}=p_{kl}^{\beta}$.

\medskip \noindent
If the two proper secants of the rank one locus:
$[d_i^{\alpha},d_j^{\alpha}]$ and $[d_i^{\beta},d_j^{\beta}]$
through this marked point {\em were distinct}, we would have a
Veronese conic $v(P_1)$
in the span $[q_i,q_j,b^{\alpha},b^{\beta}]=P_3$, and consequently
a decomposition of the restriction of ${\cal R}^3_8$ to this span
into a double-plane $P_2\supset v(P_1)$ and the cone from the
other marked point of $[q_i,q_j]$ over the Veronese conic.

\medskip \noindent
From the version of Corollary 5.3 for the $(b^{\alpha})_{\alpha}$
tetrahedron, we must have two rank two points on each edge,
and $[b^{\alpha},b^{\beta}]$, which has one rank two point
in $P_2\supset v(P_1)$, must have the other on $[q_i,q_j]$
at $p_{kl}^{\gamma}$. But this is a contradiction, because
it reduces the span  $[q_i,q_j,b^{\alpha},b^{\beta}]$ to a plane
$P_2\supset v(P_1)$.

\medskip \noindent
The source of the contradiction was to suppose that the two secants
$[d_i^{\alpha},d_j^{\alpha}]$ and $[d_i^{\beta},d_j^{\beta}]$ were
distinct. Therefore, we must have: \ $d_i^{\alpha}=d_i^{\beta}$ and
$d_j^{\alpha}=d_j^{\beta}$.

\medskip \noindent
Again, by considering the restriction of ${\cal R}^3_8$ to the
plane $[q_i,q_j,b^{\alpha},b^{\beta}]$, we find the double-line
$[d_i^{\alpha},d_j^{\alpha}]$ and the two lines from $p_{kl}^{\gamma}$
to the rank one points. Thus, $[b^{\alpha},b^{\beta}]$ has no
choice but to pass through  $p_{kl}^{\gamma}$ in order to acquire
its second rank two point.

\medskip \noindent
It is clear now that no other edge but $[b^{\gamma},b^{\delta}]$
can (and does) pass through $p_{kl}^{\alpha}$. At the same time
$d_i^{\gamma}=d_i^{\delta}$ and $d_j^{\gamma}=d_j^{\delta}$,
with their pencil running through $p_{kl}^{\gamma}$. \ \ \ $\Box$

\medskip \noindent
The proof has shown more:

\begin{cor}
The two tetrahedra $(q_i)_i$ and $(b^{\alpha})_{\alpha}$ have the
same set of twelve rank two points, distributed by two on
each system of edges.
\end{cor}

\begin{lemma}
The collection of rank one points $(d_i^{\alpha});\ i,\alpha=
1,...,4$ has exactly four distinct points, to be denoted $d_s, s=1,...,4$.
Each edge $[d_s,d_t]$ has two of the above twelve marked points.

\medskip \noindent
Thus, we have three tetrahedra with the same span:

$$ P_3=[q_1,...,q_4]=[b^1,...,b^4]=[d_1,...,d_4] $$

\noindent
and each marked point lies on some trio of edges, one from each
tetrahedron.
\end{lemma}

\medskip \noindent
{\em Proof:}\ We have seen in the proof of the previous lemma that
the intersection $[q_i,q_j]\cap [b^{\alpha},b^{\beta}]=p_{kl}^{\alpha}=
p_{kl}^{\beta}$ implies the intersection  $[q_i,q_j]\cap
[b^{\gamma},b^{\delta}]=p_{kl}^{\gamma}=p_{kl}^{\delta}$. At the same
time, the rank one points labelled $d_i^{\nu}, d_j^{\nu}; \nu=1,...,4$
are just four distinct points.

\medskip \noindent
But we also have the implication $[q_k,q_l]\cap [b^{\alpha},b^{\beta}]
=p_{ij}^{\alpha}=p_{ij}^{\beta}$, and exactly four distinct points
amongst $d_r^{\alpha}, d_r^{\beta}; r=1,...,4$.

\medskip \noindent
We want to prove that the above two tetrads of rank one points are
one and the same.

\medskip \noindent
Through $p_{kl}^{\alpha}$ runs the secant $[d_i^{\gamma},d_j^{\gamma}]$
and also the secant $[d_k^{\alpha},d_k^{\beta}]$. Should they be distinct,
we would have a Veronese conic $v(P_1)$ in their span, and thereby a
plane $v(P_1)\subset P_2\subset {\cal R}^2_6$ in $[q_1,...,q_4]$. But
this is known from previous considerations to lead to a contradiction.

\medskip \noindent
A similar argument works for the other pair. \ \ \ $\Box$

\medskip \noindent
In conclusion, we have three tetrahedra in $P_3$, with edges meeting
by three in twelve points, and with faces meeting by three in sixteen
lines. This is the Reye configuration.

\medskip \noindent
Returning to the three-space where our quadrics are surfaces,
we can first match $d_i$ with $x_i^2$ by some projective transformation,
and then make three of the marked points/quadrics match marked
points/quadrics of the standard double-four, by the action of the
torus subgroup which preserves the coordinate tetrahedron. Since
Reye configurations must match this way, so must our two other
tetrahedra match the two other of the standard double-four. By
a final switch, if necessary, the theorem is proven. \ \ \ $\Box$

\medskip \noindent
{\bf Remark:}\ It may be observed, and it will appear with even more
emphasis in the next section on tangents and Grassmannians, that
what is of the essence in the double-four example is the presence
of {\em one} tetrahedron with {\em rank four} vertices $(q_i)_i$, and
with precisely {\em two rank two
points on every edge} $[q_i,q_j]$. Assuming no face meets the
rank two locus in a conic, this leads to
${\cal R}^3_8$ restricting on each face $[q_i.q_j,q_k]$ to
four lines, with a total of sixteen lines containing the twelve
rank two points. One obtains a $(12_4,16_3)$ configuration.

\medskip \noindent
In closing this section, we illustrate the fact that the {\em smoothness assumption} in the theorem is
important, by presenting a {\em double-five} example. The example has {\em all quadrics singular, with a common
singularity} and thus actually belongs to the space of {\em conics} $P_5=P(Sym_C(3))$. It will be described as
such.

\medskip \noindent
{\bf A double-five configuration: conics.}\ {\em One can find
two quintets of conics $(q_i)_i$ and $(b_i)_i, \ i=0,...,4$,
such that any line $[q_i,b_j]$ meets the rank one locus.}

\medskip \noindent
{\em Construction:}\ Our two quintets will span the same three-space
$P_3\subset P_5=P(Sym_C(3))$. This three-space should contain a
Veronese conic $v(P_1)\subset P_2$ and just one other rank one point.

\medskip \noindent
For specificity, we'll choose the projective subspace $P_3$ of
$3\times 3$ symmetric matrices $S$ defined by:\ $s_{13}=s_{23}=0$.
Our Veronese conic is then given by:

$$ v(P_1)=\{ S\ :\ s_{11}s_{22}=s_{12}^2, \ s_{33}=0\}\subset P_2=
\{ S\ :\ s_{33}=0\}\subset P_3 $$

\noindent
and the only other point of rank one is $S_{\infty}:\ s_{33}=1$, rest 0.

\medskip \noindent
We denote it so, because it will lie on the plane at infinity with respect to an  {\em affine piece} we are
about to consider. First, we write the Veronese conic as:

$$ \frac{1}{2}(s_{11}+s_{22})^2=s_{12}^2+\frac{1}{2}(s_{11}-s_{22})^2 $$

\medskip \noindent
and then define the {\em affine piece} by $\frac{1}{2}(s_{11}+s_{22})=1$.

\medskip \noindent
In order to be closer to Euclidean intuition, we change coordinates to:

$$ x=\frac{1}{\sqrt{2}}(s_{11}-s_{22}), \ y=s_{12}, \ z=s_{33} $$

\noindent
so that the Veronese conic becomes the unit circle in the plane $z=0$.

\medskip \noindent
Since we need (and use) only {\em real points} for our configuration,
there should be no confusion here if we express the two real coordinates
$(x,y)$ by a complex number. This facilitates indicating our choice
of three points on the unit circle as the roots of unity of order three:

$$ \omega_k=e^{\frac{2\pi i}{3}k}, \  k=1,2,3; \ \ \mbox{with}\ \
\sum_{k=1}^3 \omega_k=0 $$

\medskip \noindent
Now, we can present our quintets:

$$ q_0=(0,-1),\  q_k=(2(\omega_i+\omega_j),-1),\ k=1,2,3; \
q_4=(0,\frac{1}{3})\in C\times R=R^3\subset P_3 $$

$$ b_0=(0,1),\  b_k=(2(\omega_i+\omega_j),1),\ k=1,2,3; \
b_4=(0,-\frac{1}{3})\in C\times R=R^3\subset P_3 $$

\medskip \noindent
It is elementary to verify that all lines $[q_i,b_j]$ meet the rank one
locus. \ \ \ $\Box$

\section{Tangents and Grassmannians}

\medskip \noindent
In this section we look at tangents to quadric surfaces as points in the
Grassmann variety of lines in $P_3$, that is, two dimensional vector
subspaces in $C^4$:\ $G(2,4)\subset P_5=P(\wedge^2(C^4))$. This translates
questions about common tangents into questions about intersections of
quadrics in $P_5$. In particular, we find that the variety of common
tangents for two smooth quadric surfaces in general position is
a K3 surface with 16 nodes in $G(2,4)$, more precisely:
a {\em Kummer surface}.

\medskip \noindent
We begin in arbitrary dimension:\ $P_{n-1}=P(C^n)$. Again, using the
standard bilinear form on $C^n$, we identify quadrics $q$ in $P_{n-1}$
with symmetric matrices $Q\in P(Sym_C(n))=P_{{n\choose 2}-1}$.

\medskip \noindent
The Grassmann variety $G(k,n)$ of all $(k-1)$-projective subspaces
in $P_{n-1}$, i.e. of all $k$-dimensional vector subspaces of $C^n$,
can be realized in the projective space of the $k^{th}$ {\em exterior
power} $P(\wedge^k(C^n))=P_{{n\choose k}-1}$ as all points corresponding
to {\em decomposable $k$-vectors}, that is:

$$ G(k,n)=\{ x\in P(\wedge^k(C^n))\ :\ x=\lambda x_1\wedge ...\wedge x_k
\ \mbox{for some independent set} \ x_i\in C^n\} $$

\noindent
Obviously, a $k$-subspace in $C^n$ is represented by
$x_1\wedge ...\wedge x_k$, for any choice of basis $(x_i)_i$, since
a change of basis merely introduces a proportionality factor given
by the determinant of the transition matrix.

\medskip \noindent
The conditions expressing the fact that an exterior vector, which is,
in general, a linear combination of decomposable vectors, has actually a
decomposable form are called the {\em Grassmann-Pl\"{u}cker relations},
and are all quadratic. The above realization is also called the
{\em Pl\"{u}cker embedding} of the Grassmannian $G(k,n)$. \cite{GH}

\medskip \noindent
In general: \ $dim_C G(k,n)=k(n-k)$; thus $G(2,4)\subset P_5$ is a
smooth quadric fourfold.

\medskip \noindent
{\bf Definition:}\ Let $q$ be a quadric in $P_{n-1}$. A $(k-1)$-projective
subspace $P_{k-1}\subset P_{n-1}$, is called {\bf tangent} to $q$ (at
$x\in q\cap P_{k-1}$), when the restriction of $q$ to $P_{k-1}$ i.e.
$q\cap P_{k-1}$ is singular (at $x$).

\medskip \noindent
With $q$ seen as a symmetric operator $Q$ on $C^n$, we can define a symmetric operator $\nu_k(Q)=\wedge^kQ$ on
$\wedge^k(C^n)$ by:

$$ \nu_k(Q)(x_1\wedge ...\wedge x_k)=Qx_1\wedge ...\wedge Qx_k $$

\noindent
In other words, $\nu_k(Q)$ is a quadric $\nu_k(q)$ in $P(\wedge^k(C^n))$.

\begin{lemma}
A projective subspace $P_{k-1}\subset P_{n-1}$ is tangent to the quadric $q$ in $P_{n-1}$ if and only if the
corresponding point of the Grassmannian $G(k,n)\subset P(\wedge^k(C^n))$ lies on the quadric $\nu_k(q)$.
\end{lemma}

\medskip \noindent
{\em Proof:}\ The induced standard bilinear form on $\wedge^k(C^n)$ is:

$$ <x_1\wedge ...\wedge x_k, y_1\wedge ... \wedge y_k>=det(<x_i,y_j>)_{ij} $$

\noindent
For $(x_i)_i$ a basis in our $C^k$ with $P_{k-1}=P(C^k)$, we have:

$$ <x_1\wedge ...\wedge x_k,\nu_k(Q)(x_1\wedge ...\wedge x_k)>=
det(<x_i,Qx_j>)_{ij} $$

\noindent
But the matrix $(<x_i,Qx_j>)_{ij}$ is precisely the restriction of
$q$ to our $P_{k-1}$, expressed in the chosen basis. The lemma follows.
\ \ \ $\Box$

\begin{cor}
The variety of $(k-1)$-projective subspaces tangent to a quadric
$q$ in $P_{n-1}$ is the quadratic section of the Grassmannian
$G(k,n)\subset P(\wedge^k(C^n))$ given by $\nu_k(q)$. \ \ \ $\Box$
\end{cor}

\medskip \noindent
{\bf Remark:}\ As we shall see in more detail for $\nu=\nu_2$, the map
$\nu_k$ is a {\em projection} of the $k^{th}$ Veronese map $v_k$
defined on the space $P(Sym_C(n))$ by the complete linear system of
degree $k$ hypersurfaces. In fact, $\nu_k$ corresponds to the linear
subsystem of all $k\times k$ minors, with base locus made
of quadrics of rank less than $k$.

\medskip \noindent
{\em We now fix $k=2$ and $n=4$, and thereby return to quadric surfaces.}

\medskip \noindent
We let $e_i,\ i=1,...,4$ denote the standard basis in $C^4$, and
$e_{ij}=e_i\wedge e_j,\ i<j$ the associated standard basis in
$\wedge^2(C^4)=C^6$.

\medskip \noindent
The condition for an exterior 2-vector:\ $x=\sum \alpha_{ij}e_i\wedge e_j$
to have decomposable form reads:\ $x\wedge x =0\in \wedge^4(C^4)=C$, and
in our standard basis $e_{ij}$ gives the {\em quadric}:

$$ g=2(x_{12}x_{34}-x_{13}x_{24}+x_{14}x_{23}) $$

\noindent and in matrix form $G$, with $G^2=I_6$. This is the {\em Grassmann-Pl\"{u}cker quadric}, with:

$$ G(2,4)=\{x\in P_5=P(\wedge^2(C^4))\ :\ g(x)=<x,Gx>=0\} $$

\medskip \noindent
The {\em rational map}:\ $\nu=\nu_2: P_9= P(Sym_C(4)) \cdot \cdot
\rightarrow P_{20}=P(Sym_C(6))$ takes a symmetric $4\times 4$ matrix $Q$
to the symmetric $6\times 6$ matrix $\nu(Q)$ with entries made of all
$2\times 2$ minors of $Q$. Clearly, at the projective level $\nu$ is
only defined away from the rank one locus $R^1_3\subset P_9$, and one
can eliminate the indeterminacy by {\em blowing-up} this locus.

\medskip \noindent
Since the components of $\nu$ are quadratic, it can also be presented
as a {\em quadratic Veronese map} $v_2: P_9\rightarrow P_{{9+2\choose 2}-1}=
P_{54}$ followed by some projection.

\medskip \noindent
The direction (or center) of
this projection shall be the linear span of the image $v_2(R^1_3)$.
We may recall that $R^1_3$ is itself the image of a quadratic Veronese
map $v=v_2: P_3\rightarrow P_9$, hence $v_2(R^1_3)=v_2(v_2(P_3))=v_4(P_3)$
is the image of the {\em quartic Veronese map} on $P_3$, which spans a
projective subspace of $P_{54}$ of dimension ${3+4\choose 3}-1=34$.

\medskip \noindent
Thus, the projection actually takes place on a $P_{19}$, indicating
the fact that the image of $\nu$ lies in a hyperplane of $P_{20}=
P(Sym_C(6))$. The ensuing set-up is described in:

\begin{prop}
There's a commuting diagram of regular and rational maps:

$$
\begin{array}{rcl}
P_3 & \stackrel{v_4}{\rightarrow} & P_{54}=P(Sym_C(10)) \\
v_2\downarrow & \stackrel{v_2}{\nearrow} & \downarrow \pi \\
P(Sym_C(4))=P_9 & \stackrel{\nu}{\rightarrow} &
P_{19}=\{ S\ :\ Tr(SG)=0\}\subset P_{20}=P(Sym_C(6))
\end{array}
$$

\noindent
where $G$ is the Grassmann-Pl\"{u}cker quadric, and $\pi$ is the
projection along the span of $v_4(P_3)$, which is a subspace
$P_{34}\subset P_{54}$.

\medskip \noindent
By blowing-up $P_9$ along the rank one locus $v_2(P_3)={\cal R}^1_3$
to $\tilde{P}_9$, and $P_{54}$ along $P_{34}=span[v_4(P_3)]$ to
$\tilde{P}_{54}$, this yields a diagram of regular maps:

$$
\begin{array}{rcccccc}
E_8 & \subset & \tilde{P}_9 & \stackrel{\tilde{v}_2}{\rightarrow} &
\tilde{P}_{54} & \stackrel{\pi}{\rightarrow} & P_{19} \\
\beta \downarrow & \   & \downarrow & \  & \downarrow & \  & \  \\
P_3 & \stackrel{v_2}{\rightarrow} & P_9 & \stackrel{v_2}{\rightarrow} &
P_{54} & \  & \
\end{array}
$$

\noindent
where $E_8$ denotes the exceptional divisor over the rank one locus,
with a $P_5$-bundle structure $\beta: E_8\rightarrow P_3$.

\medskip \noindent
Thus, $\tilde{\nu}=\pi\circ \tilde{v}_2 : \tilde{P}_9 \rightarrow
P_{19}$ produces a lifting of indeterminacies for $\nu$.
\end{prop}

\medskip \noindent
The fact that the image of $\nu$ lies in the hyperplane of $P_{20}$
defined by $Tr(SG)=0$ is obvious for diagonal quadrics $Q$, which
have diagonal $S=\nu(Q)$, and follows in general by the action
of $SO_C(4)$ which fixes the Grassmann-Pl\"{u}cker quadric $G$.
\ \ \ $\Box$

\medskip \noindent
{\bf Remark:}\ If we denote by $e,H$ and $h$, the {\em divisor classes}
defined on $\tilde{P}_9$ by $E_8$, (the pull-back by $\tilde{\nu}$ of)
a hyperplane in $P_{19}$, respectively (the pull-back of)
a hyperplane in $P_9$, we obtain the relation:

$$  H=2h-e \ \ \ \mbox{in the Picard group}\ \
Pic(\tilde{P}_9)=H^2(\tilde{P}_9,Z) $$

\medskip \noindent
It follows that a pencil in $P_9$ which meets the rank one locus
transversly in a single point, lifts to $\tilde{P}_9$ and meets
$E_8$ in a single point, and then maps by $\tilde{\nu}$ to a
{\em pencil} in $P_{19}$.

\medskip \noindent
On the other hand, a secant or a tangent to ${\cal R}^1_3\subset P_9$
lifts to a curve which is contracted to a point by $\tilde{\nu}$.
\ \ \ $\Box$

\medskip \noindent
For a {\em smooth quadric surface} $q\in P_9$, we have $q=Q_2\subset P_3$,
and we let $T(q)=T(Q_2)$ denote its tangent bundle. There's a natural
commutative diagram:

$$
\begin{array}{rcr}
P(T(q)) & \subset & Q_2\times G(2,4) \\
c \downarrow \ \  & \   & \downarrow \ \ \\
\tau(q)=\nu(q)\cap G(2,4) & \subset & G(2,4)
\end{array}
$$

\noindent
where $c$ is a {\em birational contraction} from the projectivised tangent
bundle of $Q_2$ onto the variety $\tau(q)$ of {\em lines in $P_3$ tangent
to $Q_2$}. The image of the exceptional locus consists of the two disjoint
conics which represent the two rulings of $Q_2$ in the Grassmannian.
$\tau(q)$ is {\em singular} along these two conics and transversal
codimension one sections acquire {\em nodes} when crossing them.

\medskip \noindent
Thus, one can look upon $\tau(q)$ as a {\em degenerate quadratic
line complex} \cite{Jes}.

\medskip \noindent
Now, we can turn our attention to the {\em variety of common tangents}
for two or more quadrics. We begin with a pair of (distinct) quadrics
$q_1$ and $q_2$, other than double-planes, and define:

$$ K(q_1,q_2)=G(2,4)\cap \nu(q_1)\cap \nu(q_2)\subset G(2,4)\subset P_5 $$

\medskip \noindent
This is a {\em complex surface} in the Grassmannian, made of all common
tangents for the two quadrics. When considered with its possible multiple
structure, it has degree $2^3=8$.

\medskip \noindent
We shall explore its structure in some relevant cases. We refer to
\cite{BPVdV} and \cite{GH} for background on compact complex surfaces.

\begin{theorem}
Suppose $(q_1,q_2)$ is a pair of smooth quadrics which is
{\em generic} in one of the subvarieties defined by the
following conditions:

\medskip
(i)\ \ \ the pencil $[q_1,q_2]$ meets the rank one locus;

\medskip
(ii)\ \ the pencil $[q_1,q_2]$ meets the rank two locus in two points;

\medskip
(iii)\ the pencil $[q_1,q_2]$ meets the rank two locus (in one point);

\medskip
(iv)\ the pencil $[q_1,q_2]$ is generic.

\medskip \noindent
Then, correspondingly, the surface of common tangents $K(q_1,q_2)$
has the following structure:

\medskip
(i)\ \ \ $P_1\times P_1$ embedded in $P_5$ by a complete linear
system of type $(1,2)$, and with multiplicity two;

\medskip
(ii)\ \ two irreducible components, each isomorphic with a nodal complete
intersection of two quadrics in $P_4$; the two components meet along
a skew quadrilateral and have nodes precisely at the four vertices of
this quadrilateral;

\medskip
(iii)\ a surface birational to a $P_1$-bundle over an elliptic
curve;

\medskip
(iv)\ \  a K3 surface with 16 nodes, more precisely: a Kummer surface.
\end{theorem}

\medskip \noindent
{\em Proof:}\  (i)\  The two quadrics
are one a basket for the other, and the genericity assumption
means in particular that they meet along a double conic $2C=q_1\cap
q_2$, with $C\approx P_1$.

\medskip \noindent
We observe first that if $x_1\in q_1$ is away from the common
conic $C$, then the tangent plane $T_{x_1}(q_1)$ meets $q_2$
along a smooth conic, and there are exactly two tangents from
$x_1$ to this conic, namely the two lines through $x_1$ in the
two rulings of $q_1$. Thus, these tangents are accounted for when
we consider all common tangents through points of $C$. The latter
make obviously the projectivized tangent bundle of $q_1$ (or $q_2$)
restricted to $C$:\ $P(T(q_1))_{|C}=P(T(q_2))_{|C}$, which is, in
fact, a trivial $P_1$-bundle over $C\approx P_1$.

\medskip \noindent
Thus, at the {\em reduced level} $K(q_1,q_2)_{red}=P_1\times P_1$, and
clearly one family of $P_1$'s in this product is plunged in $G(2,4)$
as a family of lines.

\medskip \noindent
To obtain that the embedding is actually of type $(1,2)$, we may look now
from the point of view of a {\em `basket sweep'} of $K(q_1,q_2)$, namely:
we consider the pencil $[q_1,q_2]=P_1$ as a parameter space of
common baskets $b$ (and limits thereof), and as we move $b\in [q_1,q_2]$,
the two rulings of $b$ (except at the rank three and rank one points of the
pencil) `sweep' $K(q_1,q_2)$ by pairs of conics in $G(2,4)$. At the singular
points of the pencil we have a single rational curve of common tangents.

\medskip \noindent
Thus, $K(q_1,q_2)$ appears as a $P_1$-bundle over a {\em double covering} of the `basket line' $[q_1,q_2]$, {\em
ramified} over the two singular quadrics (i.e. `degenerate baskets'). By irreducibility, the double covering is
itself a rational curve, and we have a $P_1$-bundle over $P_1$ presentation of $K(q_1,q_2)$, with the fibers
clearly plunged as conics in $G(2,4)$. The bundle is trivial by the identification of all fibers with the common
conic $C$.

\medskip \noindent
It follows that, with the proper ordering of the factors, the embedding
of $P_1\times P_1$ is of type $(1,2)$, and this yields degree four.
Hence, the surface $K(q_1,q_2)$ is actually the image  of this
embedding with multiplicity two.

\medskip
(ii)\ The fact that $K(q_1,q_2)$ is reducible whenever the pencil
$[q_1,q_2]$ has {\em two} rank two points is a consequence of the fact
that the pencil $[\nu(q_1),\nu(q_2)]$ contains a rank two quadric.

\medskip \noindent
Indeed, $\nu$ as a projection of a quadratic Veronese map, takes the
pencil $[q_1,q_2]$ to a {\em conic} in $P(Sym_C(6))$. This conic has two
rank one points corresponding to the rank two points on $[q_1,q_2]$.
The line through these two rank one points meets the pencil
$[\nu(q_1),\nu(q_2)]$ in a rank two point.

\medskip \noindent
The intersection $G(2,4)\cap \nu(q_1)\cap \nu(q_2)$ can therefore be presented as an intersection $G(2,4)\cap
\nu(q_1)\cap (P_4^+\cup P_4^-)$. For the generic case in this class, each component $G(2,4)\cap \nu(q_1)\cap
P_4^{\pm}$ is singular at the four points defined by the four lines in $q_1\cap q_2$. Their common part is a
skew quadrilateral with edges connecting the four singularities whenever they are not from the same ruling on
$q_1$ or, equivalently, on $q_2$.

\medskip \noindent
The presence of two components in $K(q_1,q_2)$ is also transparent
from the `basket sweep' approach, since there are two rational
curves of common baskets for $q_1$ and $q_2$ in this case.

\medskip
(iii)\ We can use again a `basket sweep'. From section 2
and the genericity assumption, we know that there's
{\em a smooth conic of common baskets}, say $B\approx P_1$.
It has two rank one points,
corresponding to the double-planes supported respectively by
each of the two planes of the rank two point in $[q_1,q_2]$.
It will have two other points of rank three. These four points
on $B$ are the `degenerate baskets'.

\medskip \noindent
For a proper basket $b\in B$, its two rulings provide two disjoint conics
on $K(q_1,q_2)$, while over the four `degenerate baskets' we'll have
a single rational curve. In fact, over the rank one points we have
precisely the tangents along the smooth conic component of $q_1\cap q_2$
which lies in the respective (double)-plane.

\medskip \noindent
Thus, $K(q_1,q_2)$ is {\em birationally equivalent} to a $P_1$-bundle
over a double covering of $B$ ramified over the four `degenerate baskets'.
The {\em irreducibility} of this double covering follows from a limit
argument with the rank two point moving towards a rank one
point and a case (i) situation. (The two rank one points on $B$ then
move towards the single rank one point in the limit, and the two
rank three points will coalesce into the one rank three point in the
limit). The double covering is therefore an {\em elliptic curve}.

\medskip
(iv)\ In this general case, the two quadrics meet along a degree four elliptic curve $E=q_1\cap q_2$. Also, by
the genericity assumption, the curve of pairs $E^*=\{ (x_1,x_2)\in q_1\times q_2\ :\ T_{x_1}(q_1)=
T_{x_2}(q_2)\}$ is an elliptic curve which projects on each factor as a smooth quadratic section $E^*_i\subset
q_i$.

\medskip \noindent
We'll elaborate on the role of {\em duality} in the next section, but
we should remark at this point that $E^*$ is simply the intersection
of the dual quadrics.

\medskip \noindent
It is convenient now to consider a {\em modification} $\tilde{K}(q_1,q_2)$
of our surface of common tangents $K(q_1,q_2)$ (which will turn out in fact
to be a resolution of singularities), by taking into account the points
of tangency:

$$ \tilde{K}(q_1,q_2)=\{ (x_1,x_2,t)\in q_1\times q_2\times K(q_1,q_2)\ :
\ x_i\in t\cap q_i, \ i=1,2\} $$

\medskip \noindent
The projection $\rho : \tilde{K}(q_1,q_2)\rightarrow K(q_1,q_2)$ is clearly an isomorphism away from the points
$t\in K(q_1,q_2)$ which are lines in one ruling of one quadric, and tangent somewhere to the other quadric.
Their number is easily counted as follows.

\medskip \noindent
Lines in one ruling of, say $q_1$, define a conic in $G(2,4)$; those
which are tangent to $q_2$ correspond to the intersection of this conic
with the quadric $\nu(q_2)$, and are four in number; all in all, there
are {\em sixteen points $t\in K(q_1,q_2)$ which are replaced by $P_1$'s
in} $\tilde{K}(q_1,q_2)$.

\medskip \noindent
Now we can look at one of the projections $\rho_i: \tilde{K}(q_1,q_2)
\rightarrow q_i$. Away from $E\cup E^*_i$, there are two points in
a fiber $\rho_i^{-1}(x)$, namely, the two tangents from $x$ to
the smooth conic $T_x(q_i)\cap q_j$. Over points in $E\cup E^*_i-E\cap
E^*_i$ there will be a single point.

\medskip \noindent
Let us consider finally one of the eight points $E\cap E^*_i$, say $y_i$.
Then the tangent plane $T_{y_i}(q_i)$ must be tangent to $q_j$ at
some different point $z_j\in q_j$, with $[y_i,z_j]$ the common tangent.
But with two points already in $q_j$, the whole line $[y_i,z_j]$ must
be in $q_j$. Thus the fiber $\rho_i^{-1}(y_i)$ is a $P_1$, and coincides
with one of the fibers of $\rho$.

\medskip \noindent
This means that $\tilde{K}(q_1,q_2)$ (via the Stein factorization of $\rho_i$)
is a resolution of the eight nodes of the double covering of $q_i$
ramified over $E\cup E^*_i$. Thus, $\tilde{K}(q_1,q_2)$ is a smooth
K3 surface.

\medskip \noindent
Using Nikulin's theorem in \cite{Nik}, we may conclude that the surface of
common tangents $K(q_1,q_2)$, obtained by contracting sixteen disjoint
rational curves on $\tilde{K}(q_1,q_2)$ to nodes, will be a
{\em Kummer surface}. In fact, we need not rely on this result of Nikulin,
because we may verify explicitly that the sum of the
sixteen exceptional divisors is divisible by two in
$Pic(\tilde{K}(q_1,q_2))=H^2(\tilde{K}(q_1,q_2),Z)$.

\medskip \noindent
Indeed, this follows from a calculation in the Picard lattice
of our K3 surface. Let us denote by $\tau_i$ the pull-back by $\rho_i$
of the hyperplane class of $P_3$, by $\epsilon_k^i;\  k=1,...,8$,
the classes of the curves contracted by $\rho_i$, and by $\gamma$ the
class of the elliptic curve given by all tangents to $E=q_1\cap q_2$.

\medskip \noindent
We've seen above that {\em together}, the curves contracted by
$\rho_1$ and $\rho_2$ amount exactly to the sixteen curves contracted
by $\rho$. Besides, we have:

$$ 2\tau_i=2\gamma + \sum_{k=1}^8 \epsilon^i_k \ \ \ \ i=1,2$$

\noindent
hence:

$$ \sum_{i=1}^2\sum_{k=1}^8 \epsilon^i_k=2(\tau_1+\tau_2-2\gamma) $$

\noindent
completing the argument. \ \ \ $\Box$

\medskip \noindent
{\bf Remark:}\ The divisibility by two condition verified above comes
from the representation of {\em Kummer surfaces as quotients of Abelian
surfaces by the involution} $z\mapsto -z$. The sixteen nodes then
correspond with the sixteen order two points on the Abelian surface.
Proceeding in the other direction, the Abelian surface is obtained
from $K(q_1,q_2)$ by considering the double covering of
$\tilde{K}(q_1,q_2)$ ramified over the sixteen exceptional curves.
On the covering, these rational curves become $(-1)$-curves and
can be contracted to smooth points. The resulting surface is Abelian.

\medskip \noindent
One can establish further relations in the {\em algebraic
lattice} of the K3 surface $\tilde{K}(q_1,q_2)$ obtained in the general
case.

\medskip \noindent
Let us denote by $\sigma$ the pull-back by $\rho$ of
the hyperplane class of $G(2,4)$. The degree $\sigma\gamma=8$ can be
found by an application of the Riemann-Hurwitz formula, and then we have:

$$ \sigma^2=8,\ \gamma^2=0,
\ (\epsilon^i_k)^2=-2,\ \sigma\epsilon^i_k=0,\ \gamma\epsilon^i_k=1 $$

\noindent
This leads to:

$$ 2\sigma=2\gamma +\sum_{i=1}^2\sum_{k=1}^{8} \epsilon^i_k $$

$$ \sigma=\tau_1+\tau_2-\gamma $$

\medskip \noindent
The elliptic curve $\gamma$ is part of an {\em elliptic fibration} of
our surface, where we find as another fiber the elliptic curve $\gamma^*$,
the proper transform of $E^*_i$ by $\rho_i$,
(and one and the same for $i=1,2$). \ \ \ $\Box$

\begin{prop}
Let $q_i,\ i=1,2,3$ be three distinct quadrics in a {\em generic pencil} $\ell\subset P_9=P(Sym_C(4))$. Then
their curve of common tangents:

$$ C(q_1,q_2,q_3)=G(2,4)\cap \nu(q_1)\cap \nu(q_2)\cap \nu(q_3) $$

\noindent
is given by the elliptic curve $\gamma$, made of tangents to the common
intersection $E=q_i\cap q_j,\ \{i,j\}\subset \{1,2,3\}$, and counted with
multiplicity two.
\end{prop}

\medskip \noindent
{\em Proof:}\ The image of our pencil by $\nu$  is a conic $\nu(\ell) \subset P(Sym_C(6))$ and any three
distinct points on it span its plane.

\medskip \noindent
Thus, regarding $q_1$ and $q_2$, we may consider ourselves in the generic case $(iv)$ of the preceding Theorem,
and move $q_3$ as we please along the rest of the pencil. Clearly $\gamma$, the curve of tangents to the
elliptic curve $E=q_1\cap q_2$ is always part of the intersection $K(q_1,q_2)\cap \nu(q_3)$. But any other
common tangent to $q_1$ and $q_2$ can be avoided as a tangent by some choice of $q_3\in [q_1,q_2]= \ell$.

\medskip \noindent
Thus, at the {\em reduced level} $C(q_1,q_2,q_3)_{red}=\gamma$. But
$C(q_1,q_2,q_3)$ has degree $2^4=16$, and therefore $\gamma$, with
$deg(\gamma)=\sigma\gamma=8$, has to be taken with multiplicity two.
\ \ \ $\Box$

\begin{cor}
Any four distinct quadrics in a pencil $\ell\subset P_9-{\cal R}^2_6$ have a continuum of common tangents. In
the generic case, the reduced locus is the elliptic curve $\gamma$.
\end{cor}

\medskip \noindent
Before we close this section, we may have a second look at the
{\em double-four} example in section 5. For either tetrad $(q_i)_i$,
or $(b_i)_i$, we have a degree sixteen curve of common tangents, made
of {\em eight conic components}, two for each common basket.
We can see an alternative reason for this abundant splitting  in the
fact that each pencil $[q_i,q_j]$, respectively $[b_i,b_j]$, meets the
rank two locus in two points.

\section{Duality}

\medskip \noindent
In this section we make explicit the {\em role of duality}.
For the general notion we refer to \cite{GKZ}, but here we need it
only in the case of quadrics.

\medskip \noindent
As in the previous section, we begin in arbitrary dimension $P_{n-1}=
P(V_n)$, where $V_n$ is a complex vector space of dimension $n$.

\medskip \noindent
The {\em dual projective space} $P^*_{n-1}$ is the space of all
{\em hyperplanes} in $P_{n-1}$, that is:\ $P^*_{n-1}=P(V^*_n)=G(n-1,V_n)$,
where $V^*_n$ is the dual vector space of $V_n$, and the Grassmannian
notation serves a context where no specific basis has been given.
We have a {\em canonical isomorphism}, also called {\em `orthogonality'
isomorphism}:

$$ G(k,V_n)=G(n-k,V^*_n),\ \ \mbox{taking} \ V_k\subset V_n \ \mbox{to}\
V_k^{\perp}=V^*_{n-k}\subset V^*_n $$

\noindent
where $V_k^{\perp}=V^*_{n-k}$ stands for the $(n-k)$-subspace of $V^*_n$
consisting of all linear functionals vanishing on $V_k$.

\medskip \noindent
In other words, any linear subspace $P_{k-1}\subset P_{n-1}$ has its
dual counterpart $P^*_{n-k-1}\subset P^*_{n-1}$:\ all hyperplanes
containing $P_{k-1}$.

\medskip \noindent
In this `base-free' context, we have self-dual linear operators:
$Q:V_n\rightarrow V^*_n,\ Q=Q^*$, rather than symmetric $n\times n$
matrices, and the correspondence with quadrics is given by:

$$ Q\mapsto q=\{ x\in P(V_n)\ :\ <x,Qx>=0\} $$

\noindent
where $<,>$ is the duality pairing of $V_n$ and $V^*_n$.

\medskip \noindent
{\bf Definition:}\ The {\bf dual} of a {\em smooth} quadric $q\subset
P_{n-1}$ is the subvariety $q^*\subset P^*_{n-1}$ consisting of all
hyperplanes tangent to $q$.

\begin{lemma}
$q^*$ is the quadric of $P^*_{n-1}=P(V^*_n)$ corresponding to
the self-dual operator $Q^{-1}:V^*_n\rightarrow V_n$.
\end{lemma}

\medskip \noindent
{\em Proof:}\ The hyperplane tangent to $q$ at $x\in q$ is $x^*=Qx$,
and it satisfies the equation $<Q^{-1}x^*,x^*>=<x,Qx>=0$.\ \ \ $\Box$

\medskip \noindent
If we denote by $Sym(V_n)$ the self-dual operators from $V_n$ to its
dual $V^*_n$, we obtain a {\em duality transformation}:

$$ P_{{n+1\choose 2}-1}=P(Sym(V_n))\cdot\cdot\stackrel{\ast}{\rightarrow}
P(Sym(V^*_n)) $$

\noindent as a {\em rational map} from the space of quadrics in $P_{n-1}=P(V_n)$ to the space of quadrics in
$P^*_{n-1}=P(V^*_n)$. Clearly $(q^*)^*=q$, so that duality is a {\em birational equivalence}, with inverse given
by duality applied on the target space.

\medskip \noindent
The relevance of duality for our concerns comes from the following, nearly
tautological fact:

\begin{prop}
Let $q$ be a smooth quadric in $P_{n-1}$, with dual $q^*\subset P^*_{n-1}$.
A projective $(k-1)$-subspace $P_{k-1}\subset P_{n-1}$ is tangent to $q$
if and only if its `orthogonal' $P^{\perp}_{k-1}=P^*_{n-k-1}\subset P^*_{n-1}$
is tangent to $q^*$.
\end{prop}

\medskip \noindent
{\em Proof;}\ $P_{k-1}$ is tangent to $q$ if and only if there's an
$x\in P_{k-1}\cap q$ with $P_{k-1}\subset T_x(q)=Qx$. The latter
condition reads: \ $y=Qx\in P^{\perp}_{k-1}\cap q^*$. \ \ \ $\Box$

\begin{cor}
Let $(q_i)_i$ be a collection of smooth quadrics in $P_{n-1}=P(V_n)$.
The subvariety of $G(k,V_n)$ consisting of their common tangent
$(k-1)$-planes is naturally identified, via the `orthogonality'
isomorphism, with the subvariety of $G(n-k,V^*_n)$ consisting of
common tangent $(n-k-1)$-planes for the collection of dual quadrics
$(q^*_i)_i$.
\end{cor}

\medskip \noindent
It will be convenient to have a {\em lifting of indeterminacies} for the birational equivalence determined by
duality on quadrics. At this point, we {\em choose a basis} in $V_n$, and use the standard bilinear form on
$C^n$ for the identification:\ $V_n=V^*_n=C^n$, and its consequent identifications:\
$P_{n-1}=P(V_n)=P(V^*_n)=P^*_{n-1}$, and \ $P(Sym(V_n))=P(Sym(V^*_n))=P(Sym_C(n))=P_N$, where $N={{n+1\choose
2}-1}$. Thus, the {\em duality transformation} becomes a {\em birational involution} \
$D:P_N\cdot\cdot\rightarrow P_N$, with $D(Q)=Q^{-1}$ on smooth quadrics.

\begin{prop}
Let ${\cal I}_N\subset P(Sym_C(n))\times P(Sym_C(n))=(P_N)^2$ be the
projective subvariety defined by:

$$ {\cal I}_N=\{ (A,B)\ :\ AB=\frac{1}{n}Tr(AB)\cdot I_n \} $$

\noindent
where $I_n$ denotes the identity $n\times n$ matrix.

\medskip \noindent
Then, the two projections $\pi_i:{\cal I}_N\rightarrow P_N$
are birational morphisms, and the duality transformation $D$ lifts to the
{\em involution} of ${\cal I}_N$ given by $D(A,B)=(B,A)$.
\end{prop}

\medskip \noindent
{\em Proof:}\ Clearly, for $(A,B)\in {\cal I}_N$ with $Tr(AB)\neq 0$,
we must have $B=A^{-1}$, and the closure of the graph of $D$ lies in
${\cal I}_N$.

\medskip \noindent
When $Tr(AB)=0$, we must have $AB=BA=0$, or equivalently: $Im(B)\subset Ker(A)$, or $Im(A)\subset Ker(B)$. Thus,
the fiber $\pi_1^{-1}(A)$ can be identified with quadrics in $P(C^n/Im(A))$, and similarly for the other
projection.

\medskip \noindent
It follows that, ${\cal I}_N$ can be presented as
a union of a dense open set $U_N=\{ (Q,Q^{-1})\ :\ Q \ \mbox{of rank}\ n\}$,
and closed subvarieties ${\cal R}^{i,j},\ i+j\leq n$:

$$ {\cal I}_N=U_N\cup \bigcup_{i+j\leq n} {\cal R}^{i,j} $$

\noindent
where:

$$ {\cal R}^{i,j}=\{ (A,B)\in {\cal I}_N\  : Im(A)\subset Ker(B),
\ rk(A)\leq i, \ rk(B)\leq j \} $$

\medskip \noindent
The lifting of indeterminacies on ${\cal I}_N$ is now plain, with
the additional relation:\ $D({\cal R}^{i,j})={\cal R}^{j,i}$.
\ \ \ $\Box$

\medskip \noindent
{\bf Remark:}\ One can verify that the singularity locus of ${\cal I}_N$
is:\ $Sing({\cal I}_N)=\bigcup_{i+j\leq n-3} {\cal R}^{i,j}$. In particular,
for $n=4,\ {\cal I}_9$ is smooth.

\medskip \noindent
{\em We now fix $n=4$, and return to the specifics of duality for
quadric surfaces.}

\medskip \noindent
Another expression of Proposition 7.2 for tangents to quadrics in $P_3$
is:

\begin{prop}
Let $Q\in P(Sym_C(4))$ be a rank four quadric in $P_3$.
Then:

$$ \nu(Q^{-1})=det(Q)\cdot G\nu(Q)G $$

\noindent
where, as in section 6, $G$ is the $6\times 6$ matrix corresponding
to the Grassmann-Pl\"{u}cker quadric $G(2,4)\subset P_5$ in the standard
basis.
\end{prop}

\medskip \noindent
{\em Proof:}\ With $\wedge^4(C^4)=C$, we have:

$$ <x_1\wedge x_2, G(y_1\wedge y_2)>=x_1\wedge x_2\wedge y_1\wedge y_2 $$

$$ x_1\wedge x_2\wedge Q^{-1}y_1\wedge Q^{-1}y_2= det(Q)y_1\wedge y_2
\wedge Qx_1\wedge Qx_2 $$

\medskip \noindent
and the statement is the rendering of the last equation in the standard
basis $e_{ij}=e_i\wedge e_j,\ i<j$. \ \ \ $\Box$

\medskip \noindent
{\bf Remark:}\ Conjugation with $G$ is {\em not} part of the $PSL_C(4)$
action on quadrics in $P_5$, just as the orthogonality isomorphism
$\perp: G(2,4)\rightarrow G(2,4)$ is {\em not} part of the action of
that group on $G(2,4)$.

\medskip \noindent
It may be observed that duality transforms a {\em generic pencil} in $P_9$
into a {\em rational normal cubic} which meets the rank one locus in
four points. Such relations reflect relations in the cohomology of
${\cal I}_9$.

\begin{prop}
Let $H_i$ denote the pull-back of the hyperplane class  by the
modification morphism $\pi_i:{\cal I}_9\rightarrow P_9$. Then:

$$ 4H_1={\cal R}^{3,1}+2{\cal R}^{2,2}+3{\cal R}^{1,3} $$

$$ 4H_2={\cal R}^{1,3}+2{\cal R}^{2,2}+3{\cal R}^{3,1} $$

\noindent
with the consequence:

$$ H_1+H_2={\cal R}^{3,1}+{\cal R}^{2,2}+{\cal R}^{1,3} $$

\end{prop}

\medskip \noindent
In particular, pencils which contain smooth quadrics and meet the rank one
locus will dualize to likewise pencils, an expression of the
{\em duality invariance of the `basket property'} relating two smooth
quadrics. \ \ \ $\Box$

\section{Common tangents to four spheres in $R^3$}

\medskip \noindent
In this section we interpret our results on common baskets for
the particular case of the family $P_4^{(s)}\subset P_9=P(Sym_C(4))$
which contains all quadrics whose {\em real points} are {\em spheres in}
$R^3$. Then, considering only {\em real tangents} to spheres, we determine
all {\em degenerate configurations of four spheres in $R^3$, that is:
configurations with infinitely many common tangents.}

\medskip \noindent
The generic case of configurations with finitely many common tangents,
has been studied in \cite{MPT} and \cite{ST}. The effective upper
bound 12 is in fact the complex count, which we review in the sequel.

\medskip \noindent
The affine equation of a sphere in $R^3$ is:

$$ \sum_{i=1}^3 (x_i-x_i^0)^2=r^2 $$

\noindent
with $c=(x_1^0,x_2^0,x_3^0)\in R^3$ the center, and $r=|r|>0$ the radius.

\medskip \noindent
This gives in $P_3$ the quadric:

$$ Q=\left( \begin{array}{cccc}
a_0 & 0 & 0 & a_1 \\
0 & a_0 & 0 & a_2 \\
0 & 0 & a_0 & a_3 \\
a_1 & a_2 & a_3 & a_4
\end{array} \right)
$$

\noindent
with $c=-\frac{1}{a_0}(a_1,a_2,a_3)$ and
$r^2=\frac{1}{a_0^2}(a_1^2+a_2^2+a_3^2-a_0a_4), \ \ \ a_0\neq 0$.

\medskip \noindent
We are thus led to the {\em complex projective subspace}:
$P_4^{(s)}\subset P_9=P(Sym_C(4))$ consisting of all quadrics
of the above form $Q$, with $a=(a_0:...a_4)\in P_4$. For $a_0\neq 0$,
we shall continue to designate the expressions given for $c$ and
$r^2$ as the ``center and squared radius'' of $Q$.

\medskip \noindent
Using $(x_1:x_2:x_3:x_4)$ as homogeneous coordinates in $P_3$, the
family $P_4^{(s)}$ can also be described as the family of all
quadrics in $P_3$ passing through the ``imaginary conic at infinity'':

$$ x_4=0, \ \ \  <x,x>= x_1^2+x_2^2+x_3^2=0 $$

\begin{lemma}
The subgroup of transformations in $PGL_R(4)$ preserving the family
$P_4^{(s)}$ and its real structure consists of all similarities of
$R^3=\{ x=(x_1:x_2:x_3:1)\in P_3 \}$ with respect to the above
inner product $<\ ,\ >$, that is: compositions of isometries and
rescalings. Its complexification consists of all transformations
in $PGL_C(4)$ which take the ``imaginary conic at infinity'' to
itself. \ \ \ $\Box$
\end{lemma}

\begin{lemma}
The formula:\ $ det(Q)=-a_0^2(a_1^2+a_2^2+a_3^2-a_0a_4) $ shows that
the locus $P_4^{(s)}\cap {\cal R}_8^3$ decomposes into a quadric and
the double-hyperplane $a_0^2=0$.

\medskip \noindent
In $P_9$, the three-space $a_0=0$ is tangent to the rank one locus
at the point given by $a=(0:0:0:0:1)$, which is the only rank one
point in $P_4^{(s)}$. All rank two points in $P_4^{(s)}$ are on
$a_0=0$. \ \ \ $\Box$
\end{lemma}

\begin{prop}
Let $q_i,\ i=1,2,3$ be three distinct quadrics of rank at least three
in $P_4^{(s)}$. Suppose there's a common basket for $q_i,\ i=1,2,3$.
Then:

(i)\ \ the span $[q_1,q_2,q_3]$ contains the rank one point
$T=(0:0:0:0:1)$;

\noindent
or, equivalently:

(ii)\ the centers $c_i$ of $q_i,\ i=1,2,3$ are collinear.

\medskip \noindent
A generic triple satisfying these conditions has a common basket.
\end{prop}

\medskip \noindent
{\em Proof:}\ We have to interpret Proposition 3.3 for $P_4^{(s)}$.

\medskip \noindent
Considering that the rank two locus in $P_4^{(s)}$
is $a_0=0$, our triple must satisfy the generic condition in $(C^3_{19})$
for a tangent, or condition $(H^3_{14})$. Both imply (i).

\medskip \noindent
The equivalence of (i) and (ii) is a simple computation, and the generic
converse follows from section 3. \ \ \ $\Box$

\begin{prop}
Let $q_i,\ i=1,2,3,4$ be four distinct quadrics of rank at least three
in $P_4^{(s)}$. Suppose there's a common basket for $q_i,\ i=1,2,3,4$.

\medskip \noindent
Then, the four quadrics $q_i$ lie on a conic tangent to the rank two
locus $a_0=0$ at the rank one point $T=(0:0:0:0:1)\in P_4^{(s)}$,
and the four centers are collinear.

\medskip \noindent
A generic quadruple $(q_i)$ satisfying this property has a common basket.
\end{prop}

\medskip \noindent
{\em Proof:}\ We have to interpret Proposition 4.6 for $P_4^{(s)}$.

\medskip \noindent
We see that we must be in case $(G^4_{20})$ for a tangent,
or $(J^4_{15})$. This yields the condition, with the conic degenerating to
a double line through $T$ in case $(J^4_{15})$. The fact that the four
centers must be collinear is obvious from the previous result on triples.

\medskip \noindent
The generic converse is covered by constructions in section 4. \ \ \ $\Box$

\medskip \noindent
{\bf Remark:}\ For the generic case above, we have $[q_1,...,q_4]=P_2$.
Using the triangle $T,q_1,q_2$ as simplex of reference in this plane, and
with:

$$ q_3=\alpha_0 T+\alpha_1 q_1+\alpha_2 q_2 \ \ \ \mbox{and} \ \ \
   q_4=\beta_0 T +\beta_1 q_1 +\beta_2 q_2  $$

\noindent
the existence of the conic amounts to:

$$ \alpha_0(\frac{1}{\alpha_1}+\frac{1}{\alpha_2})=
\beta_0(\frac{1}{\beta_1}+\frac{1}{\beta_2}) $$

\medskip \noindent
{\em We turn now to the problem of understanding the variety of common tangents to four spheres in $R^3$}.

\medskip \noindent
At the {\em complex projective level}, the corresponding four complex
quadrics in $P_3$ have a common curve ``at infinity'' i.e. in $x_4=0$,
namely the ``imaginary conic'': $<x,x>=0$. Tangents to this conic are
common tangents and define a conic in the Grassmannian $G(2,4)$. Thus,
what has to be identified is the remaining part of the variety of common
(complex) tangents.

\medskip \noindent
At the {\em real level}, we have to consider {\em only the real points
of this residual complex piece}, because there's no real tangent at
infinity.

\medskip \noindent
Thus, one is led to coordinates particularly adapted to lines in
the affine part $R^3\subset P_3(R)$, respectively $C^3\subset P_3$.

\medskip \noindent
A line $\ell$ in $R^3$ is completely characterized by the pair $(p,v)\in R^3\times P(R^3)$, where $p$ is the
orthogonal projection of the origin in $R^3$ on the given line, and $v$ the projective point determined by a
direction vector along the line. (One may represent $P(R^3)$ as the plane at infinity for $R^3$, and then $v$ is
simply the point where the (completed) line meets infinity.) Clearly:

$$ <p,v>=\sum_{k=1}^3 p_kv_k=0 $$

\medskip \noindent
Over $C$, the same description works generically. The resulting relation
with $G(2,4)$ is expressed in:

\begin{prop}
There's a natural birational equivalence:

$$ G(2,4)\approx I_4\subset P_3\times P_2, \ \ \  \ell\mapsto (p,v) $$

\noindent
where $I_4$ is the $P_2$-bundle over $P_2$ defined by:

$$ I_4=\{ (p,v)\in P_3\times P_2\ | \ <p,v>=\sum_{k=1}^3 p_kv_k=0 \} $$

\medskip \noindent
Let $\Gamma_4$ denote the closed graph of this birational map:

$$ \Gamma_4=\{ (\ell,p,v)\in G(2,4)\times P_3\times P_2\ | \
p\in \ell, v\in \ell, <p,v>=0 \} $$

\noindent (Here $v\in \ell$ is to be understood via the identification of directions with points at infinity:\
$P_3=C^3\cup P_2$.)

\medskip \noindent
The projection $\Gamma_4 \rightarrow G(2,4)$ is a modification
over lines at infinity (i.e. in $x_4=0$) and
lines through the origin of $C^3\subset P_3$
with null direction (i.e. $<v,v>=0$).

\medskip \noindent
The projection $\Gamma_4 \rightarrow I_4$ is a blow-up of the
rational curve $\{ (p,v)\in I_4\ | \ p_4=0,\ p=v \ \mbox{as points
at infinity}\ \}$.
\end{prop}

\medskip \noindent
{\bf Remark:}\ The fibers of $\Gamma_4 \rightarrow G(2,4)$ over
tangents to the imaginary conic at infinity are unions of two
rational curves with a common point, while elsewhere one-dimensinal
fibers are rational curves. This eventually relates to the contribution
of this conic in $G(2,4)$ in counting the isolated common tangents
to four spheres by other techniques (cf. \cite{Ful}).

\medskip \noindent
For our approach, the relevant fact in the above set-up is that
the composition $G(2,4)\cdot \cdot \rightarrow I_4 \rightarrow P_2$
is induced by a {\em linear projection} $P_5\cdot \cdot \rightarrow P_2$,
and lifting to $\Gamma_4$ resolves the indeterminacies of the map to
$I_4$, and hence to $P_2$ as well. \ \ \ $\Box$

\medskip \noindent
{\em We consider now four real quadrics of rank at least three, and
belonging to the family} $P_4^{(s)}$.
`Centers' and `squared radii' maintain a formal sense and, after a
translation, we may assume the centers are at $0,c_1,c_2,c_3\in R^3$,
with corresponding squared radii $r^2,r_1^2,r_2^2,r_3^2$.

\medskip \noindent
A way to set aside the component given by tangents `at infinity',
is to write the equations for common tangents in $(p,v)$ coordinates,
with $p\in C^3$. As in \cite{ST}, the equations are:

$$ <p,v>=0 $$

$$ <p,p>=r^2 $$

$$ <c_i,p>=\frac{1}{2<v,v>}[-<c_i,v>^2+<v,v>(<c_i,c_i>+<p,p>-r_i^2)] $$

\medskip \noindent

\begin{prop} Suppose the four centers are affinely independent (i.e.
the real span of $c_i,\ i=1,2,3$ is $R^3$). Then, counting multiplicity,
there are twelve complex common tangents `away from infinity' for the
four quadrics.
\end{prop}

\medskip \noindent
{\em Proof:}\ With centers understood as column vectors, we put
$M=[c_1\ c_2\ c_3]^t$. It is a {\em real invertible matrix} and
the last three equations take the form:

\medskip \noindent

$$ M p =\frac{1}{2<v,v>}[\Phi_2(v)+<v,v>\Phi_0] $$

\noindent
where $\Phi_2(v)$ is the column vector with entries $-<c_i,v>^2$, and
$\Phi_0$ is the column vector with entries $<c_i,c_i> +r^2-r_i^2$.

\medskip \noindent
Thus $v\in P_2$ determines $p$, and must satisfy:

$$ <M^{-1}(\Phi_2(v)+<v,v>\Phi_0),v>=0 $$

$$ <M^{-1}(\Phi_2(v)+<v,v>\Phi_0),M^{-1}(\Phi_2(v)+<v,v>\Phi_0)>=
4r^2<v,v>^2 $$

\medskip \noindent
We prove that there can be no one-dimensional
component in the intersection of the above cubic and quartic curves
by showing that the further intersection with the conic
$<v,v>=0$ is empty. Indeed, the equations yield the system:

$$ <v,v>=0 $$

$$ <M^{-1}\Phi_2(v),v>=0 $$

$$ <M^{-1}\Phi_2(v),M^{-1}\Phi_2(v)>=0 $$

\medskip \noindent
The first two equations say that $M^{-1}\Phi_2(v)$ is on the tangent
at $v$ to the smooth conic $<v,v>=0$, and the last that $M^{-1}\Phi_2(v)$
is itself on the same conic. This means:

$$ M^{-1}\Phi_2(v)=\mu v \ \ \ \mbox{that is:}\ \ \ \Phi_2(v)=\mu Mv $$

\medskip \noindent
But this gives:

$$ (<c_1,v>^2:<c_2,v>^2:<c_3,v>^2)=(<c_1,v>:<c_2,v>:<c_3,v>) $$

\noindent
which has only {\em real solutions}, namely: $(1:1:1)$ or $(1:1:0)$ or
$(1:0:0)$, up to permutation. In all cases, with $M$ real, the resulting $v$
is {\em real} and we cannot have $<v,v>=0$.

\medskip \noindent
The cubic and quartic curves have therefore zero-dimensional
intersection, that is, counting multiplicity, they meet in twelve
points. The twelve solutions determine twelve common tangents `away from
infinity'. \ \ \ $\Box$

\begin{cor}
Four spheres in $R^3$ with affinely independent centers have at
most twelve common real tangents.
\end{cor}

\medskip \noindent
{\bf Remark:}\  Configurations of four spheres with twelve distinct
real common tangents are constructed in \cite{MPT}. See also \cite{ST}.

\medskip \noindent
The next case to consider is when {\em the four centers are coplanar
but no three of them are collinear}. It requires more detailed
computations for ruling out the possibility of infinitely many common
tangents in the real case.

\begin{prop}
Four quadrics of rank at least three from $P_4^{(s)}(R)$, with coplanar
centers but no three of them collinear, have only isolated common tangents
`away from infinity'. Their number is at most twelve.
\end{prop}

\medskip \noindent
{\em Proof:}\ By Lemma 8.1, we may assume that the centers
span $x_3=0$. Thus:

$$ M=\left( \begin{array}{ccc}
c_{11} & c_{12} & 0 \\
c_{21} & c_{22} & 0 \\
c_{31} & c_{32} & 0
\end{array} \right)
$$

\noindent
We let $M_{12}$ stand for the $2\times 2$ upper left corner.

\medskip \noindent
Eliminating $p$ from the equations yields in this case {\em a sextic and
a conic} in $v\in P_2$, and our aim is to show that their intersection
has to be zero-dimensional. The {\em conic} $E_2$ is obtained by using a
non-zero vector $k$ in the kernel of $M^t$:

$$ \sum_{i=1}^3 k_ic_i=0, \ \ \ \ k\neq 0 $$

$$ 0 =2<v,v><Mp,k>=<\Phi_2(v),k>+<v,v><\Phi_0,k> $$

\noindent
The sextic is obtained by writing $p=p_{12}+p_3e_3$, with $p_{12}\perp e_3$.
Then, with similar notations for $v=(v_{12},v_3),\ \Phi_2(v)$ etc. :

$$ p_3=-\frac{<p_{12},v_{12}>}{v_3} \ \ \ \mbox{from} \ \ <p,v>=0 $$

$$ p_{12}=\frac{1}{2<v,v>}M_{12}^{-1}[\Phi_2(v)_{12} +<v,v>\Phi_{0,12}] $$

\medskip \noindent
With $\Psi_2(v)=\Psi_2(v_{12})=M_{12}^{-1}\Phi_2(v)_{12}$,\
$\Psi_0=M_{12}^{-1}\Phi_{0,12}$, and $||x||^2=<x,x>\in C$ this gives the
{\em sextic} $E_6$:

$$ v_3^2||\Psi_2(v)+<v,v>\Psi_0||^2 +
<\Psi_2(v)+<v,v>\Psi_0,v_{12}>^2-4<v,v>^2v_3^2r^2=0 $$

\medskip \noindent
Considering that our centers $c_i,\ i=1,2,3$ are in the plane of the
first two coordinates, we shall envisage them as two-dimensional vectors
when this simplifies formulae. Thus, from:

$$ \frac{1}{k_3}k_{12}=-(M_{12}^t)^{-1}c_3 $$

\noindent
we obtain an equivalent expression for the conic $E_2$:

$$ <\Psi_2(v),c_3>+<c_3,v_{12}>^2-<v,v><\Phi_0,k>=0 $$

\medskip \noindent
From here on, our proof relies on various computational consequences
of the above equations for the sextic $E_6$ and conic $E_2$, which we
present in a sequence of lemmas.

\begin{lemma}
$E_6,E_2$ and $<v,v>=0$ have no common solution $v\in P_2$, unless $v_{12}=c_i^{\perp}=(c_{i2}:-c_{i1})$ as
points in $P_1$, for some $i\in \{1,2,3\}$, and $<c_i^{\perp},c_j-c_k>=0$. (The last condition means that the
four centers are the vertices of a trapeze.)
\end{lemma}

\medskip \noindent
{\em Proof:}\ With $<v,v>=0$, $E_6$ and $E_2$ become equations in
$v_{12}\in P_1$:

$$ ||v_{12}||^2||\Psi_2(v_{12})||^2-<\Psi_2(v_{12}),v_{12}>^2=0 \ \ \ \
\ (E_6^{12}) $$

$$ <\Psi_2(v_{12}),c_3>+<c_3,v_{12}>^2=0 \ \ \ \ \ \ \ (E_2^{12}) $$

\medskip \noindent
The first equation requires:

$$ \Psi_2(v_{12})=\lambda v_{12} \ \ \ \ \  i.e.\ \ \ \ \ \Phi_2(v)_{12}=
\lambda M_{12}v_{12} $$

\noindent
which gives with $E_2^{12}$:

$$ <c_i,v>^2=-\lambda <c_i,v> \ \ \ \mbox{for} \ \ i=1,2,3 $$

\noindent
Any two centers being independent, $v_{12}$ must be orthogonal to some
$c_i$ and then $<c_i^{\perp},c_j-c_k>=0$. \ \ \ $\Box$

\begin{lemma}
For $E_2$ and $E_6$ to have a common one-dimensional component,
it is necessary that $E_2^{12}$ and $E_6^{12}$ have two common
solutions.
\end{lemma}

\medskip \noindent
{\em Proof:}\ Suppose there's a single solution $v_{12}=c_i^{\perp}$.
By relabelling, if necessary, we may assume $i=3$.
Then the common component of $E_6$ and $E_2$
must be the line in $P_2$ through $(c_3^{\perp}:\pm i<c_3,c_3>^{1/2})$,
or, for $<c_3,c_3>=0$, the tangent $<c_3,v>=0$ to $<v,v>=0$. In either case
it's the line through $(c_3^{\perp}:0)$ and $(0:0:1)$.

\medskip \noindent
This means that, when we rewrite $E_6$ as an equation in $v_3$ (actually
$v_3^2$), with coefficients depending on $v_{12}$, all these coefficients
must vanish identically for $v_{12}=c_3^{\perp}$. In other words, we put
$<v,v>=<v_{12},v_{12}>+v_3^2$ in $E_6$ and obtain a cubic in $v_3^2$.
The vanishing of its four coefficients for $v_{12}=c_3^{\perp}$ yields:

$$ ||\Psi_0||^2=4r^2 $$

$$ 2<\Psi_2(c_3^{\perp}),\Psi_0>  +<\Psi_0,c_3^{\perp}>^2=0 $$

$$ ||\Psi_2(c_3^{\perp})||^2 + 2||c_3||^2<\Psi_2(c_3^{\perp}),\Psi_0>=0 $$

$$ <\Psi_2(c_3^{\perp}),c_3^{\perp}> + ||c_3||^2<\Psi_0,c_3^{\perp}>=0 $$

\medskip \noindent
From the previous lemma we know that:

$$ \Psi_2(c_3^{\perp})= \mu c_3^{\perp} \ \ \  \mbox{with} \ \ \
\mu = -<c_1,c_3^{\perp}> = -<c_2,c_3^{\perp}> $$

\noindent
and the last three equations become:

$$ <\Psi_0,c_3^{\perp}>(<\Psi_0,c_3^{\perp}>+2\mu)=0 $$

$$ ||c_3||^2(\mu+2<\Psi_0,c_3^{\perp}>)=0 $$

$$ ||c_3||^2(<\Psi_0,c_3^{\perp}>+\mu)=0 $$

\medskip \noindent
With real quadrics, we have $||c_3||^2\neq 0$, and the last
two equations already provide a contradiction, since $\mu \neq 0$.\ \ \
$\Box$

\medskip \noindent
Again, by relabelling the centers, we may assume that the two
common solutions of $E_2^{12}$ and $E_6^{12}$ are
$v_{12}=c_1^{\perp}$ and $v_{12}=c_2^{\perp}$. Thus, $c_1+c_2=c_3$
i.e. the centers form a parallelogram.

\medskip \noindent
This suggests using a {\em translation} which brings the origin at
the center of the parallelogram. We assume therefore that the four
centers are now at $a=(a_1,a_2,0)^t, b=(b_1,b_2,0)^t, -a,$ and
$-b$, with squared radii $r_i^2,\ i=1,...,4$.

\medskip \noindent
The original system becomes:

$$ <p,v>=0 $$

$$ <a,p>=\frac{1}{2<v,v>}[-<a,v>^2+<v,v>(<a,a>+<p,p>-r_1^2)] $$

$$ -<a,p>=\frac{1}{2<v,v>}[-<a,v>^2+<v,v>(<a,a>+<p,p>-r_3^2)] $$

$$ <b,p>=\frac{1}{2<v,v>}[-<b,v>^2+<v,v>(<b,b>+<p,p>-r_2^2)] $$

$$ -<b,p>=\frac{1}{2<v,v>}[-<b,v>^2+<v,v>(<b,b>+<p,p>-r_4^2)] $$

\medskip \noindent
Subtraction in the last two pairs of equations gives:

$$ 2<a,p>=r_3^2-r_1^2 $$

$$ 2<b,p>=r_4^2-r_2^2 $$

\medskip \noindent
This shows that the first two coordinates $p_{12}$ of $p$ are
determined by centers and squared radii alone, and remain constant. But
this means that {\em all common tangents to our four real quadrics
meet the perpendicular drawn from $p_{12}$ to the plane of the
centers.}

\medskip \noindent
{\bf Remark:} A theorem in \cite{MS} already addresses a situation of this nature, and shows that the common
tangents to three spheres which meet at the same time a fixed line cannot be infinitely many {\em unless the
three spheres have collinear centers} (and the fixed line adequate position). However, it is not necessary to
rely on this result in order to prove our proposition, as we show next.

\medskip \noindent
Completing to an equivalent system, we have:

$$ p=p_{12} - \frac{<p_{12},v_{12}>}{v_3}e_3 $$

$$<a,v>^2=<v,v>[<a,a>+<p,p>-\frac{1}{2}(r_3^2+r_1^2)] $$

$$<b,v>^2=<v,v>[<b,b>+<p,p>-\frac{1}{2}(r_4^2+r_2^2)] $$

\medskip \noindent
With:\ $\alpha=<a,a>-\frac{1}{2}(r_3^2+r_1^2)$ and
$\beta=<b,b>-\frac{1}{2}(r_4^2+r_2^2)$, we obtain:

$$ \frac{<a,v>^2}{<b,v>^2}=\frac{\alpha+<p,p>}{\beta+<p,p>} $$

\medskip \noindent
For $v\in P_2$ the system amounts now to intersecting a conic and a quartic:

$$ <a+b,v><a-b,v>=<v,v>[<a+b,a-b> +\frac{1}{2}(r_4^2+r_2^2-r_3^2-r_1^2)] $$

$$ \frac{<a,v_{12}>^2}{<b,v_{12}>^2}=
\frac{(\alpha+<p_{12},p_{12}>)v_3^2+<p_{12},v_{12}>^2}
{(\beta+<p_{12},p_{12}>)v_3^2+<p_{12},v_{12}>^2} $$

\medskip \noindent
With:\ $A=\alpha+<p_{12},p_{12}>,\ \ B=\beta+<p_{12},p_{12}>$ and
$C=A-B=\alpha-\beta$, the equations say:

$$ v_3^2=\frac{1}{C}<a+b,v><a-b,v> - <v_{12},v_{12}> $$

$$ v_3^2=\frac{<p_{12},v_{12}>^2<a+b,v_{12}><a-b,v_{12}>}
{A<b,v_{12}>^2-B<a,v_{12}>^2} $$

\medskip \noindent Thus, for the conic and quartic to have a common
 one-dimensional component it is necessary that:

$$ <a+b,v_{12}><a-b,v_{12}>[C<p_{12},v_{12}>^2 + B<a,v_{12}>^2 -
A<b,v_{12}>^2]= $$

$$ = -C<v_{12},v_{12}>[A<b,v_{12}>^2-B<a,v_{12}>^2] $$

\noindent identically in $v_{12}\in P_1$.

\medskip \noindent
Now, evaluating at $v_{12}=(a+b)^{\perp}$ and $(a-b)^{\perp}$, we find:
$A=B=C=0$

\medskip \noindent
Returning these conditions into the system gives:

$$ <a+b,v_{12}><a-b,v_{12}>=0 $$

$$ <a,v_{12}>^2=<b,v_{12}>^2=(<v_{12},v_{12}>+v_3^2)p_3^2 $$

\noindent
The first equation requires:\ $v_{12}=(a+b)^{\perp}$ or $(a-b)^{\perp}$, and
then the second determines $v_3$, since we cannot have $p_3=0$ with $a$ and $b$
linealy independent. Thus, there's no one-dimensional
family of solutions. \ \ \ $\Box$

\begin{cor}
Four spheres with coplanar centers but no three of them collinear
have at most twelve common real tangents.
\end{cor}

\medskip \noindent
Finally, when {\em three of the centers are collinear}, we have rotational
symmetry around this axis for the common tangents to the corresponding
three quadrics. Thus, either (i) the three quadrics have a common conic
in the affine part $C^3$, or (ii) the three quadrics have a common basket
(and only one by Proposition 8.3 and Lemma 3.2).

\medskip \noindent
Accordingly, the fourth quadric cannot have a curve of common tangents
with the other three in the affine part unless it passes through the
same common conic, in case (i), or has the same common basket in case
(ii).

\medskip \noindent
Both cases require the four centers to be collinear, and, restricting to
the case of spheres and real tangents, we obtain the result described in
the introduction:

\begin{theorem}
Four distinct spheres in $R^3$ have infinitely many common real tangents
if and only if they have collinear centers and at least one common real
tangent.

\medskip \noindent
This means that either all four spheres intersect in a circle, possibly
degenerating to a common tangency point, or each sphere has a curve of
tangency with one and the same real quadric of revolution with symmetry axis
determined by the line passing through all centers. This quadric can be
a cone, a cylinder, or a one-sheeted hyperboloid.
\end{theorem}

\medskip \noindent
{\bf Remark:}\ Our argument has made effective use of reality assumptions.
It will be observed that the complex case allows more possibilities for
degenerate configurations.

\medskip
Ciprian Borcea \ \ \ \ \ \ \ borcea@rider.edu

Rider University, Lawrenceville, NJ 08648, USA.

\medskip
Xavier Goaoc \ \ \ \ \ \ \ \ goaoc@loria.fr

Sylvain Lazard \ \ \ \ \ \ \ \  lazard@loria.fr

Sylvain Petitjean \ \ \ \ \ \ petitjea@loria.fr

LORIA-INRIA Lorraine, CNRS, Univ. Nancy 2, France.
\end{document}